\documentclass[psamsfonts,reqno]{amsart}
\usepackage{amssymb,eucal,graphics,xy}
\xyoption{all}
\xyoption{ps}
\numberwithin{equation}{section}

\newcommand{\pup}[1]{\textup{(}#1\textup{)}}

\theoremstyle{plain}

\newtheorem{lemma}{Lemma}[section]
\newtheorem{theorem}[lemma]{Theorem}
\newtheorem{proposition}[lemma]{Proposition}
\newtheorem{corollary}[lemma]{Corollary}
\newtheorem{examplepf}[lemma]{Example}
\newtheorem{claim}{Claim}

\newtheorem*{stat}{\name}
\newcommand{\name}{testing}

\theoremstyle{definition}
\newtheorem{definition}[lemma]{Definition}
\newtheorem{example}[lemma]{Example}
\newtheorem{problem}{Problem}

\theoremstyle{remark}
\newtheorem{remark}[lemma]{Remark}

\newenvironment{all}[1]{\renewcommand{\name}{#1}\begin{stat}}
                        {\end{stat}}

\newcommand{\qedc}{{\qed}~{\rm Claim~{\theclaim}.}}

\newenvironment{cproof}
{\begin{proof}[Proof of Claim.]}
{\qedc\renewcommand{\qed}{}\end{proof}}

\newcommand{\set}[1]{\left\{#1\right\}}
\newcommand{\setm}[2]{\set{{#1}\mid{#2}}}

\DeclareMathOperator{\supp}{supp}
\DeclareMathOperator{\Irr}{Irr}
\DeclareMathOperator{\Gen}{Gen}
\DeclareMathOperator{\Id}{Id}
\DeclareMathOperator{\Idc}{Id_c}
\DeclareMathOperator{\Conc}{Con_c}

\newcommand{\CH}{\mathcal{H}}
\newcommand{\CP}{\mathcal{P}}
\newcommand{\arch}{\mathrm{arch}}
\newcommand{\lear}{\leq^{\arch}}
\newcommand{\slar}{<^{\arch}}

\newcommand{\pe}{\propto}
\newcommand{\as}{\asymp}
\newcommand{\pear}{\pe^{\arch}}
\newcommand{\asar}{\as^{\arch}}

\newcommand{\es}{\varnothing}

\newcommand{\two}{\boldsymbol{2}}
\newcommand{\three}{\boldsymbol{3}}

\newcommand{\xe}{\boldsymbol{e}}
\newcommand{\xf}{\boldsymbol{f}}
\newcommand{\xg}{\boldsymbol{g}}
\newcommand{\xh}{\boldsymbol{h}}
\newcommand{\xs}{\boldsymbol{s}}
\newcommand{\xt}{\boldsymbol{t}}

\newcommand{\bx}{\boldsymbol{x}}
\newcommand{\by}{\boldsymbol{y}}

\newcommand{\bu}{\boldsymbol{u}}
\newcommand{\bv}{\boldsymbol{v}}
\newcommand{\bw}{\boldsymbol{w}}

\newcommand{\jz}{\ensuremath{\langle\vee,0\rangle}}
\newcommand{\jzu}{\ensuremath{\langle\vee,0,1\rangle}}
\newcommand{\jzl}{\ensuremath{\langle\vee,\wedge,0\rangle}}
\newcommand{\jzh}{\jz-ho\-mo\-mor\-phism}
\newcommand{\jzuh}{\jzu-ho\-mo\-mor\-phism}
\newcommand{\jzlh}{\jzl-ho\-mo\-mor\-phism}
\newcommand{\jzs}{\jz-se\-mi\-lattice}
\newcommand{\eps}{\varepsilon}
\newcommand{\cm}{commutative monoid}
\newcommand{\poag}{partially ordered abelian group}
\newcommand{\povs}{partially ordered vector space}
\newcommand{\pss}{pseu\-do-sim\-plic\-i\-al}
\newcommand{\Pss}{Pseu\-do-sim\-plic\-i\-al}
\newcommand{\NN}{\mathbb{N}}
\newcommand{\ZZ}{\mathbb{Z}}
\newcommand{\QQ}{\mathbb{Q}}

\newcommand{\LL}{\mathcal{L}}
\newcommand{\MM}{\mathcal{M}}

\newcommand{\SSS}{\mathcal{S}}
\newcommand{\SSfb}{\mathcal{S}_{\mathrm{fb}}}
\newcommand{\EE}{\mathcal{E}}
\newcommand{\EED}{\mathcal{E}_{\mathrm{d}}}
\newcommand{\EEps}{\mathcal{E}_{\mathrm{ps}}}

\begin{document}

\title[Lifting semilattice diagrams]%
{Liftings of diagrams of semilattices\\
by diagrams of dimension groups}

 \author[J.~T\r uma]{Ji\v r\'\i\ T\r uma}
 \address{Department of Algebra\\
          Faculty of Mathematics and Physics\\
          Sokolovsk\'a 83\\
          Charles University\\
          186 00 Praha 8\\
          Czech Republic}
 \email{tuma@karlin.mff.cuni.cz}

 \author[F.~Wehrung]{Friedrich Wehrung}
 \address{CNRS, UMR 6139\\
          Universit\'e de Caen, Campus II\\
          D\'epartement de Math\'ematiques\\
          B.P. 5186\\
          14032 CAEN Cedex\\
          FRANCE}
 \email{wehrung@math.unicaen.fr}
 \urladdr{http://www.math.unicaen.fr/\~{}wehrung}
\thanks{The authors were supported by the \textsc{Barrande} program and
by the institutional grant CEZ:J13/98:113200007a. The first author was
also partially supported by GA CR 201/99 and by GA UK 162/1999}


\maketitle

\section*{Introduction}
There are various ways to obtain distributive semilattices from other
mathematical objects. Two of them are the following; we refer to
Section~\ref{S:Basic} for more precise definitions.
A \emph{dimension group} is a directed, unperforated \poag\ with the interpolation
property, see also K.\,R. Goodearl \cite{Gpoag}. With a dimension group $G$ we
can associate its semilattice of \emph{compact} ( = finitely generated)
\emph{ideals}
$\Idc G$. Because of the interpolation property the positive cone $G^+$ has the
refinement property, thus the compact ideal semilattice $\Idc G$
is distributive. A ring $R$ is \emph{locally matricial} over a field $K$, if it
is isomorphic to a direct limit of finite products of full matricial rings
$\MM_n(K)$. If $R$ is locally matricial, or, more generally, \emph{regular} (in
von~Neumann's sense), then its semilattice $\Idc R$ of compact ideals is
distributive, see, for example, K.\,R. Goodearl \cite{Gvnrr}.

These two different contexts are related as follows. With a locally matricial
ring~$R$, we can associate its (partially ordered) \emph{Grothendieck group}
$K_0(R)$. It turns out that $K_0(R)$ is a dimension group, and the following
relation holds:
 \begin{equation}\label{Eq:IdG=IdR}
 \Idc R\cong\Idc(K_0(R)),
 \end{equation}
see \cite[Corollary~15.21]{Gvnrr}. Therefore, the compact ideal semilattice of a
locally matricial ring is the compact ideal lattice of a dimension group. By the
results of K.\,R. Goodearl and D.\,E. Handelman \cite{GoHa86}, the converse holds
for dimension groups of cardinality at most~$\aleph_1$, but it fails for larger
dimension groups, see F. Wehrung \cite{Wehr98,Wehr99}. Nevertheless this makes it
worthwhile to study the functor that with a dimension group~$G$ associates its
semilattice of compact ideals $\Idc G$. For a survey paper on this and related
issues we refer the reader to K.\,R. Goodearl and F. Wehrung \cite{GoWe}.

The central open problem that motivates the present paper is Problem~10.1
of~\cite{GoWe}, that asks whether any distributive \jzs\ of cardinality at
most~$\aleph_1$ is isomorphic to the compact ideal semilattice of some dimension
group. Let us call this the \emph{compact ideal representation problem}. By
\cite{GoHa86}, this is equivalent to asking whether any distributive
\jzs\ of cardinality at most $\aleph_1$ is isomorphic to the compact ideal
semilattice of some locally matricial ring. Although we do not solve this problem
here, we settle related issues with some interest of their own, and, hopefully,
that will provide a stepping stone towards a solution for the general problem.

We first recall the answers to some related questions.
As a corollary of the main result of G.\,M. Bergman~\cite{Berg86} and a
theorem of K.\,R. Goodearl \cite{Gvnrr}, the following result is obtained, see
also \cite{GoWe} for different proofs.

\begin{theorem}\label{T:Count}
Every countable distributive \jzs\ $S$ is isomorphic to $\Idc G$ for
some countable dimension group $G$.
\end{theorem}

Another result of \cite{GoWe} is the following.

\begin{theorem}\label{T:Latt}
Every distributive lattice $L$ with zero
is isomorphic to $\Idc G$ for some dimension group $G$.
\end{theorem}

On the other hand, P.~R\r{u}\v{z}i\v{c}ka \cite{Ruzi1} constructs a distributive
\jzs\ of cardinality $\aleph_2$ that is not isomorphic to $\Idc G$ for any
dimension group $G$.

These results leave open the case of distributive \jzs s of
cardinality~$\aleph_1$. By P. Pudl\'ak \cite{Pu}, every distributive
\jzs\ is isomorphic to the direct union of its finite distributive
\jz-subsemilattices. A variant of this result, stating that every distributive
\jzs\ is a direct limit of finite Boolean \jzs s, is proved in \cite{GoWe}. Thus
it is natural to try an inductive approach to the compact ideal representation
problem for distributive semilattices of cardinality~$\aleph_1$:
find a \emph{simultaneous representation} of a suitable
limit system of finite Boolean semilattices that converges to a given
distributive semilattice $S$ of cardinality $\aleph_1$. By A.\,P.~Huhn
\cite{Hu}, one can assume that the limit system is indexed by an infinite lattice
in which every principal ideal is  a finite \emph{dismantlable} lattice.

After preliminary results, most notably about the existence of the so-called
$\lambda$-generic maps in Sections \ref{S:FlGenSp} and \ref{S:FlGen}, we
prove in Section~\ref{S:DismLift} our main results:

\begin{all}{Theorem~\ref{T:dislift}}
Every finite dismantlable diagram of finite Boolean \jzs s has a
simultaneous representation with respect to the functor $\Idc$.
\end{all}

\begin{all}{Theorem~\ref{T:Nab1-1}}
For any countable distributive \jzs\ $S$ and any countable dimension vector space
$H$, every \jzh\ $\xf\colon S\to\Idc H$ can be lifted by a positive homomorphism
$f\colon G\to H$ for some dimension vector space $G$.
\end{all}

We also state consequences of Theorem~\ref{T:Nab1-1}, as well in ring theory as
in lattice theory, see Corollaries~\ref{C:Bab1-1ring} and \ref{C:Bab1-1latt}.

The remaining sections contain various counterexamples to other strategies for a
positive solution of the compact ideal representation problem for distributive
semilattices of cardinality $\aleph_1$.

\section{Basic concepts}\label{S:Basic}

For basic concepts about \poag s, we refer the reader to K.\,R. Goodearl
\cite{Gpoag}. We recall some of the definitions here.
For a \poag\ $G$, we denote by $G^+=\setm{x\in G}{0\leq x}$ the positive cone
of $G$, and we say that $G$ is \emph{directed}, if $G=G^++(-G^+)$. We put
$G^{++}=G^+\setminus\set{0}$, and we note $\NN=\ZZ^{++}$.
A subgroup $H$ of $G$ is \emph{convex}, if $0\leq x\leq a$
implies that $x\in H$, for all $a\in H$ and $x\in G$. We say that $H$ is an
\emph{ideal} of $G$, if it is a directed, convex subgroup of $G$, and a
\emph{compact ideal} of $G$, if it is an ideal generated by a finite subset
of $G$.

As every ideal of $G$ is directed, any ideal (resp., compact ideal) of $G$
is generated (as an ideal) by a subset of $G^+$ (resp., \emph{finite} subset
of $G^+$). In the latter case, one may replace a finite subset
$\set{a_1,\dots,a_n}$ of $G^+$ by the singleton $\set{a}$, where
$a=\sum_{i=1}^na_i$, so the compact ideals of $G$ are exactly the subsets of
$G$ of the form
 \[
 G(a)=\setm{x\in G}{\exists n\in\ZZ^+\text{ such that }-na\leq x\leq na},
 \qquad\text{for }a\in G^+.
 \]
An \emph{order-unit} of $G$ is an element $a\in G^+$ such that $G(a)=G$.
We denote by $\Id G$ the lattice of ideals of $G$, and by $\Idc G$ the
\jzs\ of compact ideals of~$G$. Observe that $\Id G$ is an algebraic lattice
and that $\Idc G$ is its \jzs\ of compact elements, see
G. Gr{\"a}tzer \cite{Grat98} for unexplained terminology. For elements $a$,
$b\in G^+$, we write $a\pe b$, if there exists $n\in\NN$ such that $a\leq nb$, and
$a\as b$, if $a\pe b\pe a$. Hence
 \[
 \Idc G=\setm{G(a)}{a\in G^+},
 \]
with containment and equality among the $G(a)$-s determined by
 \[
 G(a)\subseteq G(b)\text{ if{f} }a\pe b\qquad\text{and}\qquad
 G(a)=G(b)\text{ if{f} }a\as b,
 \]
for all $a$, $b\in G^+$. Hence, $\Idc G$ is isomorphic to $\nabla(G^+)$, the
\emph{maximal semilattice quotient of $G^+$}, see \cite{GoWe}.

The assignment $G\mapsto\Idc G$ can be naturally extended to a
\emph{functor}, by defining, for a positive homomorphism (that is, an
order-preserving group homomorphism) $f\colon G\to H$ of \poag s and a
compact ideal $I$ of $G$, $(\Idc f)(I)$ as the compact ideal of $H$
generated by the image $f[I]$ of $I$ under $f$. Hence, $(\Idc f)(G(a))=H(f(a))$,
for all $a\in G^+$, and it is an easy exercise to verify that the functor $\Idc$
thus defined \emph{preserves direct limits}.

Most of the \poag s we shall deal with will be equipped with an additional
structure of \emph{vector space}, always over the field $\QQ$ of rational
numbers. Any \povs\ $E$ has the additional property that $mx\geq0$ implies
that $x\geq0$, for all $m\in\NN$ and all $x\in E$, we say that $E$ is
\emph{unperforated}, see \cite{Gpoag}. We observe that all group
homomorphisms between vector spaces preserve the vector space structure.

A \emph{refinement monoid} (see A. Tarski \cite{Tars}, F. Wehrung \cite{Wehr92})
is a \cm\ $M$ such that for any positive integers $m$, $n$ and elements $a_0$,
\ldots, $a_{m-1}$, $b_0$, \ldots, $b_{n-1}$ of $M$ such that
$\sum_{i<m}a_i=\sum_{j<n}b_j$, there are elements $c_{i,j}$ (for $i<m$
and $j<n$) such that $a_i=\sum_{j<n}c_{i,j}$, for all $i<m$, and
$b_j=\sum_{i<m}c_{i,j}$, for all $j<n$.

A \poag\ $G$ is an \emph{interpolation group} (see \cite{Gpoag}) if
for all positive integers $m$, $n$ and elements $a_0$, \ldots,
$a_{m-1}$, $b_0$, \ldots, $b_{n-1}$ of $G$ such that $a_i\leq b_j$,
for all $i<m$ and $j<n$, there exists $x\in G$ such that
$a_i\leq x\leq b_j$, for all $i<m$ and $j<n$. Equivalently, the positive cone
$G^+$ is a refinement monoid (see \cite[Proposition 2.1]{Gpoag}). In
addition, we say that $G$ is a \emph{dimension group} (see
\cite{Gpoag}) if~$G$ is directed and unperforated. If $G$ is an
interpolation group, then $\Idc G$ is a refinement semilattice. Refinement
semilattices are usually called \emph{distributive semilattices}, see
\cite{Grat98}.

We can prove right away the following very elementary lifting result:

\begin{proposition}\label{P:LiftSG}
For every \jzs\ $S$, there exists a \povs\ $E$ such that $\Idc E\cong S$.
\end{proposition}

\begin{proof}
We first observe that $S$ can be \jz-embedded into some powerset semilattice
$\CP(X)$, for a set $X$. For example, take $X=S$ and use the embedding that with
any $s\in S$ associates $j(s)=\setm{x\in S}{s\nleq x}$. We then let $F$ denote
the vector space of all maps $f\colon X\to\QQ$ with finite range. For any
$f\in F$, we put $\supp f=\setm{x\in X}{f(x)\neq 0}$, and we put
 \[
 M=\setm{f\in F}{f(x)\geq0,\text{ for all }x\in X\text{ and }\supp f\in S}.
 \]
Then $M$ is an additive submonoid of $F$ (in particular, it is cancellative),
$M\cap(-M)=\set{0}$, and $M$ is closed under multiplication by
nonnegative scalars. Hence $M$ is the positive cone of a structure of directed
\povs\ on $E=M+(-M)$. Let $\pi\colon M\to S$ be the map defined by the rule
$\pi(f)=\supp f$, for all $f\in M$. Then $\pi(f)\leq\pi(g)$ (resp.,
$\pi(f)=\pi(g)$) if{f} $f\pe g$ (resp.,
$f\as g$) in $E$, for all $f$, $g\in M$, whence $\pi$ induces an isomorphism
$\varphi$ from $\Idc E$ onto $S$ by the rule $\varphi(E(f))=\supp f$, for all
$f\in M$.
\end{proof}

Unfortunately, even if $S$ is distributive, the \povs\ $E$ constructed above may
not have interpolation.

Throughout the paper we shall often formulate our results in basic
categorical language. Among the main categories we shall be working with are
the following:

\begin{itemize}
\item[---] The category $\EE$ of pairs $A=(A_0,1_A)$, where $A_0$ is a
\povs\ (that we will subsequently identify with $A$) and $1_A$ is an
order-unit of $A_0$. For objects $A$ and $B$ of $\EE$, a
\emph{homomorphism}, or \emph{positive homomorphism}, from $A$ to $B$ is
a positive homomorphism $f\colon A_0\to B_0$. Observe that we
do \emph{not} require that $f(1_A)=1_B$. We say that a homomorphism
$f\colon A\to B$ is \emph{normalized}, if $f(1_A)=1_B$. In particular, we
will slightly abuse terminology by using isomorphisms only in the `normalized'
sense, that is, an isomorphism $f\colon A\to B$ satisfies the equality
$f(1_A)=1_B$.

\item[---] The full subcategory $\EED$ of $\EE$ consisting of all
\emph{dimension vector spaces} in~$\EE$.

\item[---] The category $\SSS$ of \jzs s and \jzh s.

\item[---] The full subcategory $\SSfb$ of $\SSS$ whose objects are the
finite Boolean semilattices.
\end{itemize}
In particular, $\Idc$ defines, by restriction, a functor from $\EE$ to
$\SSS$ that preserves direct limits.

Every partially ordered set $(P,\leq)$ can be viewed as a category with
objects $x\in P$ and a unique morphism $\eps_{x,y}\colon x\to y$ whenever
$x\leq y$ in $P$. Let $\EE'$ be a subcategory of $\EE$. We say that a
diagram $\Phi\colon(P,\leq)\to\SSS$ of semilattices (i.e., a functor
from $P$, viewed as a category, to $\SSS$) has a
\emph{lifting}, with respect to the functor $\Idc$, to the category $\EE'$,
if there is a diagram $\Psi\colon(P,\leq)\to\EE'$ such that the two functors $\Phi$
and $\Idc\circ\Psi$ are naturally equivalent, that is, if there are
isomorphisms $\iota_x\colon \Phi(x)\to\Idc(\Psi(x))$, for $x\in P$, such that the
following diagram commutes,
 \[
{
\def\labelstyle{\displaystyle}
\xymatrix{
 \Idc(\Psi(x))\ar[rr]^{\Idc(\Psi(\eps_{x,y}))} & & \Idc(\Psi(y))\\
 \Phi(x)\ar[u]^{\iota_x}\ar[rr]_{\Phi(\eps_{x,y})} & & \Phi(y)\ar[u]_{\iota_y}
}
}
 \]
for all $x\leq y$ in $P$. We will usually omit the phrase \emph{with respect to}
$\Idc$, since lifts with respect to other functors will not be investigated.

For a set $I$, we denote by $\CP(I)$ the powerset of $I$, and, for
any $i\in I$, we put
 \[
 \CP_i(I)=\setm{X\in\CP(I)}{i\in X},\qquad
 \CP_i^*(I)=\setm{X\in\CP(I)}{i\notin X}.
 \]

\section{Pseudo-simplicial spaces}\label{S:Pseudosimp}

For a nonempty, finite set $X$, let $\QQ_X$ denote the \povs\ with
underlying vector space $\QQ^X$ and partial ordering defined by
 \[
 f\leq g\Leftrightarrow\text{ either }f=g\text{ or }
 f(x)<g(x),\text{ for all }x\in X,
 \]
endowed with the canonical order-unit $1^X$, defined as the constant
function on $X$ with value $1$. Then $\QQ_X$ is a dimension vector
space. It is \emph{simple}, that is,
$\Idc\QQ_X\cong\two$. We denote by
$(1^X_x)_{x\in X}$ the canonical basis of $\QQ^X$, so that
$1^X=\sum_{x\in X}1^X_x$. We observe that while $1^X\in(\QQ_X)^{++}$, the
vectors $1^X_x$, for $x\in X$, do not belong to $(\QQ_X)^{++}$ unless $X$
is a singleton. Nevertheless they keep a certain positivity character,
captured by the following definition:

\begin{definition}\label{D:ArchOrd}
For a \poag\ $G$, we define a partial preordering $\lear$ on $G$ by the rule
 \begin{equation}\label{Eq:Defarch}
 x\lear y\Longleftrightarrow\exists u\in G\text{ such that }
 nx\leq ny+u,\text{ for all }n\in\ZZ^+,
 \end{equation}
for all $x$, $y\in G$. Then we define the partially preordered abelian group
$G_\arch=(G,\lear)$, the \emph{archimedean quotient} of $G$.
\end{definition}

In particular, for any finite set $X$, the archimedean quotient of $\QQ_X$ is
$\QQ^X$, and $0\slar 1^X_x$, for all $x\in X$.

\begin{definition}\label{D:PseudSimp}
An object $A$ of $\EE$ is
\begin{enumerate}
\item \emph{simple \pss}, if $A\cong\QQ_X$ for some nonempty finite set $X$,

\item \emph{\pss}, if $A=\bigoplus_{i<n}A_i$ for some natural
number $n$ and some simple \pss\ spaces $A_i$, for $i<n$.
\end{enumerate}
\end{definition}

The formula in (ii) above means that $A=\bigoplus_{i<n}A_i$ as vector spaces,
$1_A=\sum_{i<n}1_{A_i}$, and the canonical bijection from $\prod_{i<n}A_i$ to
$A$ is an isomorphism in~$\EE$ (that is, for the vector space structure as well
as the---componentwise---partial order structure).

In particular, \emph{simplicial vector spaces}, that is, spaces of the
form $\QQ^n$ with the positive cone $(\QQ^n)^+$
consisting of all vectors with non-negative coordinates, 
are particular cases of \pss\ spaces. Observe that for any \pss\ space~$E$,
the space $E_{\arch}$ is simplicial.

For a simple \pss\ space $A$, the \emph{canonical basis} of $A$ is defined
as the set $T$ of minimal elements $t$ of $A_\arch^{++}$ such that
$t\wedge(1_A-t)=0$, where $\wedge$ denotes the meet operation in the
lattice-ordered group $A_\arch$.
Hence $T$ is a (finite) vector space basis of $A$ and the map from $\QQ_T$
to $A$ that with any sequence
$(x_t)_{t\in T}$ associates $\sum_{t\in T}x_tt$ is an isomorphism from
$\QQ_T$ onto $A$. It sends $1^T_t$ to $t$, for all $t\in T$, and $1^T$ to $1_A$.

We denote by $\EEps$ the full subcategory of $\EE$ whose objects are the
\pss\ spaces.

We state without proof the following lemma, that summarizes some elementary
properties of \pss\ spaces.

\begin{lemma}\label{L:DecPss}
Let $m$ be a natural number, let $A_0$, \dots, $A_{m-1}$ be simple
\pss\ spaces. Put $A=\bigoplus_{i<m}A_i$, and identify $A_i$ with its canonical
image in $A$, for all $i<m$. Then the following assertions hold.

\begin{enumerate}
\item The ideals of $A$ are exactly the subsets of the form
$\bigoplus_{i\in I}A_i$ for a subset $I$ of $\set{0,1,\dots,m-1}$.

\item The simple ideals of $A$ are exactly the $A_i$, for $i<m$.

\item $\Idc A\cong\two^m$, a finite Boolean semilattice.
\end{enumerate}
\end{lemma}

For a \povs\ $E$, the two binary relations $\pe$ and $\as$ on~$E^+$
can be refined as follows. 
For a positive, rational number $\lambda$ and $a$, $b\in E^+$, we
introduce the following notations:
 \begin{align*}
 a\pe_\lambda b,&\quad\text{if }a\leq\lambda b;\\
 a\as_\lambda b,&\quad\text{if }a\pe_\lambda b\text{ and }
 b\pe_\lambda a;\\
 a\pear_\lambda b,&\quad\text{if }a\pe_{\lambda+\eps}b,\quad
 \text{for all rational }\eps>0;\\
 a\asar_\lambda b,&\quad\text{if }a\as_{\lambda+\eps}b,\quad
 \text{for all rational }\eps>0.
 \end{align*}
Hence $a\pear_\lambda b$ (resp., $a\asar_\lambda b$) in $E$ implies that
$a\pe_\lambda b$ (resp., $a\as_\lambda b$) in $E_\arch$; the converse holds if
$E$ is \emph{simple}.

\section{Basic facts about refinement}\label{S:BasicRef}

We start with an easy and useful lemma.

\begin{lemma}\label{L:MultRef}
Let $M$ be a refinement monoid, let $I$ and $T_i$, $i\in I$, 
be finite nonempty sets, let $a$, $a_{i,j}$ \pup{for $i\in I$ and $j\in T_i$} be
elements of $M$ such that
 \[
 a=\sum_{j\in T_i}a_{i,j},\quad\text{for all }i\in I.
 \]
Then there are elements $x_\varphi$ \pup{for $\varphi\in T=\prod_{i\in I}T_i$}
of $M$ such that
 \begin{equation}\label{Eq:Decaij}
 a=\sum_{\varphi\in T}x_\varphi\quad\text{and}\quad
 a_{i,j}=\sum_{\substack{\varphi\in T\\
 \varphi(i)=j}}x_\varphi,
 \quad\text{for all }\quad i\in I\quad\text{and}\quad j\in T_i.
 \end{equation}
\end{lemma}

\begin{proof}
An easy induction on the cardinality of $I$. There is nothing to prove if
$|I|=1$. Now assume that the claim holds for any set $I$ of cardinality
$n\geq1$ and $I'=I\cup\set{k}\neq I$. By the induction hypothesis on $I$ we
get elements  $x_\varphi\in M$ for $\varphi\in T$ such that
 \[
 a=\sum_{\varphi\in T}x_\varphi\quad\text{and}\quad
 a_{i,j}=\sum_{\substack{\varphi\in T\\ \varphi(i)=j}}x_\varphi,
 \quad\text{for all}\quad i\in I\quad\text{and}\quad j\in T_i.
 \] 
Since also $a=\sum_{l\in T_k}a_{k,l}$ and $M$ is a refinement monoid, we get
that there are elements $x_{\varphi,l}\in M$, for $\varphi\in T$
and $l\in T_k$, such that
 \[
 x_{\varphi}=\sum_{l\in T_k}x_{\varphi,l}\quad\text{and}\quad
 a_{k,l}=\sum_{\varphi\in T}x_{\varphi,l}
 \] 
for any $\varphi\in T$ and $l\in T_k$. Thus
 \[
 a=\sum_{l\in T_k}a_{k,l}=\sum_{l\in T_k}\sum_{\varphi\in T}x_{\varphi,l}.
 \] 
Since $T\times T_k=\prod_{i\in I\cup\set{k}}T_i$, the indices $(\varphi,l)$ are
in one-to-one correspondence with functions $\psi\in T\times T_k$. Moreover,
for every $i\in I$ and $j\in T_i$ we get
 \[
 a_{i,j}=\sum_{\substack{\varphi\in T\\ \varphi(i)=j}}x_\varphi=
 \sum_{\substack{\varphi\in T\\ \varphi(i)=j}}\sum_{l\in T_k}x_{\varphi,l}=
 \sum_{\substack{\psi\in T\times T_k\\ \psi(i)=j}}x_\psi.
 \] 
Finally, for every $l\in T_k$ we get
 \begin{equation}
 a_{k,l}=\sum_{\varphi\in T}x_{\varphi,l}=
 \sum_{\substack{\psi\in T\times T_k\\ \psi(k)=l}}x_\psi.\tag*{\qed}
 \end{equation}
\renewcommand{\qed}{}
\end{proof}

\begin{lemma}\label{L:LamAs}
Let $\lambda\geq1$ be a rational number, let $E$ be a dimension
vector space, let $I$ be a finite set, let $a_i$ \pup{for $i\in I$}
be elements of $E^+$. Then the following are equivalent:
 \begin{enumerate}
 \item $a_i\as_\lambda a_j$, for all $i$, $j\in I$.

 \item There are elements $b_X$ \pup{for $X\in\CP(I)$} of $E^+$ such that
 \[
a_i=\sum_{X\in\CP_i^*(I)}b_X+\lambda\cdot\sum_{X\in\CP_i(I)}b_X,\quad
 \text{for all }i\in I.
 \]
 \end{enumerate}
\end{lemma}

\begin{proof}
(ii)$\Rightarrow$(i) is trivial.

Now assume that (i) holds. By applying interpolation to the system of
inequalities $\lambda^{-1}a_i\leq a_j$, for all $i$, $j\in I$, we obtain $a\in E$
such that $\lambda^{-1}a_i\leq a\leq a_j$, for all $i$,
$j\in I$. Hence, $a\leq a_i\leq\lambda a$, for all $i\in I$, so that
there are $a'_i$, $a''_i\in E^+$ such that $a_i=a+a'_i$ and
$(\lambda-1)a=a'_i+a''_i$, for all $i\in I$. By Lemma~\ref{L:MultRef}
applied to $M=E^+$ and $J=\set{0,1}$, there are elements $b_X$ of $E^+$,
for $X\in\CP(I)$, such that
 \begin{gather*}
 a'_i=(\lambda-1)\sum_{X\in\CP_i(I)}b_X,\qquad
 a''_i=(\lambda-1)\sum_{X\in\CP_i^*(I)}b_X,
 \qquad\text{for all }i\in I;\\
 a=\sum_{X\in\CP(I)}b_X.
 \end{gather*}
Therefore, for all $i\in I$,
 \begin{align*}
 a_i&=\sum_{X\in\CP(I)}b_X+(\lambda-1)\sum_{X\in\CP_i(I)}b_X\\
 &=\sum_{X\in\CP_i^*(I)}b_X+\lambda\cdot\sum_{X\in\CP_i(I)}b_X.
 \tag*{\qed}
 \end{align*}
\renewcommand{\qed}{}
\end{proof}

\section{Flat and generic homomorphisms with simple target}
\label{S:FlGenSp}

\begin{definition}\label{D:FlGen}
Let $m$ be a natural number, let $A_0$, \ldots, $A_{m-1}$, $B$ be
\emph{simple} objects of~$\EED$, let
$f\colon\bigoplus_{i<m}A_i\to B$ be a positive homomorphism, let
$\lambda\geq1$ be a rational number. We say that $f$ is
\begin{enumerate}
\item \emph{$\lambda$-flat}, if $f(1_{A_i})\asar_\lambda f(1_{A_j})$,
for all $i$, $j<m$ such that $f(1_{A_i})$, $f(1_{A_j})\neq0$.

\item \emph{$\lambda$-generic}, if for any simple \pss\
space $C$ and for any $\lambda$-flat homomorphism
$g\colon\bigoplus_{i<m}A_i\to C$ such that
 \[
 f[A_i]\neq\set{0}\Leftrightarrow g[A_i]\neq\set{0},\quad\text{for all }i<m,
 \]
there exists a homomorphism $h\colon B\to C$ such that $g=h\circ f$.
\end{enumerate}
\end{definition}

The following simple result summarizes some of the basic properties of
flat and generic homomorphisms:

\begin{lemma}\label{L:FlGen}
Let $m$ be a natural number, let $A$ be a \pss\ space,
let $A_0$, \ldots, $A_{m-1}$, $B$, $G$ be
\emph{simple} objects of~$\EED$.
Let $\alpha$, $\beta$, $\lambda$ be rational numbers such that
$1\leq\alpha\leq\beta$ and $1\leq\lambda$.
\begin{enumerate}
\item Every homomorphism $f\colon\bigoplus_{i<m}A_i\to B$ is $\lambda$-flat for
some $\lambda\geq 1$.

\item Every $\alpha$-flat homomorphism from a \pss\ space to a simple
\pss\ space is $\beta$-flat.

\item Every $\beta$-generic homomorphism from a \pss\ space to a simple
\pss\ space is $\alpha$-generic.

\item Suppose that $G$ is \pss,
let $f\colon A\to G$ and $h\colon A\to B$ be homomorphisms, and let
$\xg\colon\Idc G\to\Idc B$ be a \jzh\ such that
$\Idc h=\xg\circ\Idc f$. If $f$ is $\lambda$-generic and $h$ is
$\lambda$-flat, then there exists a homomorphism $g\colon G\to B$ such
that $\xg=\Idc g$ and $h=g\circ f$.
\end{enumerate}

\end{lemma}

The situation of (iv) above may be summarized by the following
commutative diagrams:
 \[
{
\def\labelstyle{\displaystyle}
\xymatrix{
 & \Idc B & & & B\\
 \Idc G\ar[ur]^{\xg} & & & G\ar@{-->}[ur]^g\\
 & \Idc A\ar[lu]^{\Idc f}\ar[-2,0]_{\Idc h} & & &
 A\ar[lu]^f\ar[-2,0]_h
}
}
 \]
\begin{proof}
To prove (i) observe that if $f(1_{A_i})\ne 0$, then it is an order-unit in
$B$. Thus for any $i$, $j<m$ such that both $f(1_{A_i})$ and $f(1_{A_j})$
are nonzero, there exists $\lambda_{i,j}\in\QQ^{++}$ such 
that $f(1_{A_i})\leq \lambda_{i,j}f(1_{A_j})$. Now set $\lambda$ any positive
rational number larger than~$1$ and all $\lambda_{i,j}$.

Claims (ii) and (iii) are straightforward.

To prove (iv) observe first that $\xg$ has a lifting,
that is, a homomorphism $g\colon G\to B$ such that
$\Idc g=\xg$. Indeed, if $\xg=0$, take $g=0$, while if $\xg$ is the
unique isomorphism from $\Idc G$ to $\Idc B$, any nonzero
positive homomorphism from $G$ to~$B$ is a lifting of $\xg$, for
example, if $G=\QQ_m$ and if $b\in B^{++}$, the map
$g\colon G\to B$ defined by the rule
$g(x_0,\ldots,x_{m-1})=\left(\sum_{i<m}x_i\right)b$, for all $x_0$,
\dots, $x_{m-1}\in\QQ$.

Now, for (iv) above, if $h=0$, then any lifting $g$ of $\xg$
satisfies the required conditions. So suppose that $h\neq0$. Then, by
the $\lambda$-genericity of $f$ and the $\lambda$-flatness of $h$,
there exists a homomorphism $g\colon G\to B$ such that $h=g\circ f$.
In particular, $g\neq0$. Similarly, from
$\Idc h=\xg\circ\Idc f$ and $\Idc h\neq0$ follows that
$\xg\neq0$. Therefore, since both $G$ and $B$ are simple, $g$ lifts
$\xg$.
\end{proof}

We shall now define the canonical generic maps. We are given a \pss\
space $A$ and a \jzh\ $\xf\colon\Idc A\to\two$. We shall
construct a simple \pss\ space $G=\Gen(A,\xf)$ and a lifting
$f\colon A\to G$ of $\xf$ with respect to the $\Idc$ functor.

Put $A=\bigoplus_{i<m}A_i$, for a natural number $m$ and simple
\pss\ spaces $A_0$, \ldots, $A_{m-1}$. We shall define a \pss\ space
$G=\Gen(A,\xf)$ and, for every rational number $\lambda\geq1$, a lifting
$f_\lambda\colon A\to G$ of $\xf$ with respect to the $\Idc$ functor.

For all $i<m$, let $T_i$ denote the
canonical basis of $A_i$ (see Section~\ref{S:Pseudosimp}). Hence the~$T_i$,
for $i<m$, are mutually disjoint finite sets and
$\bigcup_{i<m}T_i$ is a vector space basis of $A$. Let
$I=\setm{i<m}{\xf(A_i)=1}$ (observe that $A_i$ is a compact ideal
of~$A$, thus it belongs to the domain of
$\xf$, see Lemma~\ref{L:DecPss}). Furthermore, we put
 \[
 T=\prod_{i\in I}T_i,\qquad F=\CP(I)\times T,\qquad\text{and}\qquad
 G=\QQ_F.
 \]
Observe that both $T$ and $F$ are finite, nonempty sets.
For a rational number $\lambda\geq1$, we define a linear map
$f_\lambda\colon A\to\QQ_F$ by its action on the elements of
$\bigcup_{i<m}T_i$. For $i<m$ and $t\in T_i$, we define $f_\lambda(t)=0$ if
$i\notin I$, while, if $i\in I$, we put
 \begin{equation}\label{Eq:DefCGen}
 f_\lambda(t)=
 \sum_{\substack{(X,\varphi)\in F\\ i\notin X,\ \varphi(i)=t}}
 1^F_{(X,\varphi)}+\lambda\cdot
 \sum_{\substack{(X,\varphi)\in F\\ i\in X,\ \varphi(i)=t}}
 1^F_{(X,\varphi)}.
 \end{equation}
For all $i<m$, let $f_{i,\lambda}$ denote the restriction of $f_\lambda$ to
$A_i$.

\begin{lemma}\label{L:filPosHom}
The map $f_{i,\lambda}$ is a positive homomorphism from
$A_i$ to $\QQ_F$, for all $i<m$, and it is nonzero if{f}
$i\in I$.
\end{lemma}

\begin{proof}
For $i\notin I$, $f_{i,\lambda}=0$. Suppose now that $i\in I$.
It is trivial that $f_{i,\lambda}$ defines a positive homomorphism
from $(A_i)_\arch\cong\QQ^{T_i}$ to $\QQ^F$. To prove that it is also a
nonzero positive homomorphism from $A_i$ to $\QQ_F$, it is sufficient
to prove that $f_{i,\lambda}(1_{A_i})$ is an order-unit of $\QQ^F$. By
observing that $1_{A_i}=\sum_{t\in T_i}t$ and by using \eqref{Eq:DefCGen},
we can compute:
 \[
 f_{i,\lambda}(1_{A_i})=
 \sum_{\substack{(X,\varphi)\in F\\ i\notin X}}
 1^F_{(X,\varphi)}+\lambda\cdot
 \sum_{\substack{(X,\varphi)\in F\\ i\in X}}
 1^F_{(X,\varphi)}.
 \]
Since all the components of this vector relatively to all
$1^F_{(X,\varphi)}$, for $(X,\varphi)\in F$, are positive,
$f_{i,\lambda}(1^{T_i})$ is indeed an order-unit of $\QQ^F$.
\end{proof}

It follows from Lemma~\ref{L:filPosHom} that $f_\lambda$ is a positive
homomorphism from $A$ to $\QQ_F$.
It follows immediately from Lemma~\ref{L:filPosHom} that
$f_\lambda\colon A\to\QQ_F$
is a lifting of the semilattice
map $\xf\colon\Idc A\to\two$.
Indeed, if we take the canonical isomorphism from
$\Idc\QQ_F$ to~$\two$, the following diagram
 \[
{
\def\labelstyle{\displaystyle}
\xymatrix{
\two \ar[r]^(.4){\mathrm{can.}} & \Idc\QQ_F\\
\Idc A\ar[u]^{\xf}\ar[ru]_{\Idc f_\lambda}
}
}
 \]
is commutative, since $f_\lambda[A_i]\neq\set{0}$ if{f}
$\xf(A_i)=1$, for all $i<m$.

\begin{lemma}\label{L:flflat}
The homomorphism $f_\lambda$ is $\lambda$-flat.
\end{lemma}

\begin{proof}
For all $i\in I$, $f_\lambda(1_{A_i})=f_{i,\lambda}(1_{A_i})$ is a
linear combination of the vectors $1^F_{(X,\varphi)}$, for
$(X,\varphi)\in F$, with coefficients either $1$ or $\lambda$ (see
\eqref{Eq:DefCGen}). The conclusion follows.
\end{proof}

We observe here the relevance of the definition of flatness by using $\asar$
instead of $\as$: indeed, the vectors $1^F_{(X,\varphi)}$, for
$(X,\varphi)\in F$, do not belong to the positive cone of $\QQ_F$ except if
$F$ is a singleton.

Now we come to the main result of this section:

\begin{lemma}\label{L:flGen}
The homomorphism $f_\lambda$ is $\lambda$-generic.
\end{lemma}

\begin{proof}

Let $C$ be a simple \pss\ space, let $g\colon\bigoplus_{i<m}A_i\to C$
be a $\lambda$-flat homomorphism such that $g[A_i]\neq\set{0}$ if{f}
$f_\lambda[A_i]\neq\set{0}$ (that is, $i\in I$), for all $i<m$. Set again
$A_i=\QQ_{T_i}$ for every $i<m$. We define elements $a_i$ (for $i<m$) and
$a_{(t)}$ (for $i<m$ and $t\in T_i$) of $C$ by
 \[
 a_i=g(1_{A_i}),\qquad a_{(t)}=g(t).
 \]
Observe that $a_i\in C^+$ while we can only say that
$a_{(t)}\in C_{\arch}^+$. Moreover, $a_i=\sum_{t\in T_i}a_{(t)}$. {}From
the fact that $g$ is $\lambda$-flat follows that
$a_i\asar_\lambda a_j$ in $C$, for all $i$, $j\in I$, hence, by
Lemma~\ref{L:LamAs} applied to $C_\arch$ (that has interpolation, because
$C$ is \pss), there are elements $b_X$ (for
$X\in\CP(I)$) of $C_\arch^+$ such that
 \begin{equation}\label{Eq:Decai}
 a_i=\sum_{X\in\CP_i^*(I)}b_X+\lambda\cdot\sum_{X\in\CP_i(I)}b_X,\quad
 \text{for all }i\in I.
 \end{equation}
Now, for fixed $i\in I$, we apply refinement to the equality
 \[
 \sum_{t\in T_i}a_{(t)}=
 \sum_{X\in\CP_i^*(I)}b_X+\lambda\cdot\sum_{X\in\CP_i(I)}b_X,
 \]
which holds in $C_{\arch}^+$. 
We obtain elements $c_{X,i,t}$ (for $(X,i)\in\CP(I)\times I$ and
$t\in T_i$) of $C_{\arch}^+$ such that
 \begin{align}
 a_{(t)}&=\sum_{X\in\CP_i^*(I)}c_{X,i,t}+\lambda\cdot
 \sum_{X\in\CP_i(I)}c_{X,i,t}&&
 (\text{for all }i\in I\text{ and }t\in T_i)\label{Eq:Decait1}\\
 b_X&=\sum_{t\in T_i}c_{X,i,t}&&(\text{for all }X\in\CP(I)).
 \label{Eq:DecbX1}
 \end{align}
Now we apply Lemma~\ref{L:MultRef} to the system of equations
\eqref{Eq:DecbX1} in $C_{\arch}^+$. We find elements
$d_{(X,\varphi)}$ (for $(X,\varphi)\in F$) of $C_{\arch}^+$ such that
 \begin{align}
 b_X&=\sum_{\varphi\in T}d_{(X,\varphi)}&&
 (\text{for all }X\in\CP(I))\label{Eq:DecbX}\\
 c_{X,i,t}&=
 \sum_{\substack{\varphi\in T\\ \varphi(i)=t}}d_{(X,\varphi)}&&
 (\text{for all }X\in\CP(I),\ i\in I,\ t\in T_i).\label{Eq:DeccXit}
 \end{align}
We define a linear map $h\colon\QQ^F\to C$ by the rule
 \begin{equation}\label{Eq:Defofh}
 h(1^F_{(X,\varphi)})=d_{(X,\varphi)},
 \quad\text{for all }(X,\varphi)\in F.
 \end{equation}
Hence $h$ is a positive homomorphism from $\QQ^F$ to $C_\arch$. To
verify that $h$ is a positive homomorphism from $\QQ_F$ to $C$, it
suffices to verify that $h(1^F)\in C^+$. This is trivial for $I=\es$
(then $\QQ_F=\QQ$). Suppose that $I\neq\es$. We compute:
 \begin{align*}
 h(1^F)=\sum_{(X,\varphi)\in F}h(1^F_{(X,\varphi)})&=
 \sum_{(X,\varphi)\in F}d_{(X,\varphi)}\\
 &=\sum_{X\in\CP(I)}b_X
 \quad(\text{by }\eqref{Eq:DecbX}).
 \end{align*}
Hence, it follows from \eqref{Eq:Decai} that $h(1^F)\as a_i$ in
$C_\arch$ for any $i\in I$. But $a_i\in C^+$, thus $h(1^F)\in C^+$.

It remains to prove that $h\circ f_\lambda(t)=g(t)$,
for all $i<m$ and all $t\in T_i$. This is obvious if $i\notin I$, in
which case both sides of the equality are zero. So suppose that $i$
belongs to $I$. We compute:
 \begin{align*}
 h\circ f_\lambda(t)&=
 \sum_{\substack{(X,\varphi)\in F\\ i\notin X,\,\varphi(i)=t}}
 d_{(X,\varphi)}+\lambda\cdot
 \sum_{\substack{(X,\varphi)\in F\\ i\in X,\,\varphi(i)=t}}
 d_{(X,\varphi)}&&(\text{by }\eqref{Eq:DefCGen}\text{ and }
 \eqref{Eq:Defofh})\\
 &=\sum_{X\in\CP_i^*(I)}c_{X,i,t}+\lambda\cdot
 \sum_{X\in\CP_i(I)}c_{X,i,t}&&(\text{by }\eqref{Eq:DeccXit})\\
 &=a_{(t)}&&(\text{by }\eqref{Eq:Decait1})\\
 &=g(t)&&(\text{by the definition of }a_{(t)}),
 \end{align*}
which concludes the proof.
\end{proof}

\section{Flat and generic homomorphisms with arbitrary target}
\label{S:FlGen}

\begin{definition}\label{D:FlGen2}
Let $n$ be a natural number, let $B_0$, \dots, $B_{n-1}$ be
simple \pss\ spaces, let $A$ be a \pss\ space, let
$f\colon A\to\bigoplus_{j<n}B_j$ be a homomorphism. Write
 \[
 f(x)=(f_0(x),\ldots,f_{n-1}(x)),\qquad\text{for all }x\in A,
 \]
with $f_j\colon A\to B_j$, for all $j<n$. For any rational number
$\lambda\geq1$, we say that $f$ is
\begin{enumerate}
\item \emph{$\lambda$-flat}, if $f_0$, \dots, $f_{n-1}$ are
$\lambda$-flat.

\item \emph{$\lambda$-generic}, if $f_0$, \dots, $f_{n-1}$ are
$\lambda$-generic.
\end{enumerate}
\end{definition}

Because of Lemma~\ref{L:DecPss}, the decomposition
$B=\bigoplus_{j<n}B_i$ of $B$ into simple ideals is unique, hence the
definition above is well-stated.

Part of Lemma~\ref{L:FlGen} has a straightforward analogue for these
more general definitions of flat and generic maps:

\begin{lemma}\label{L:FlGen2}
Let $\alpha$ and $\beta$ be rational numbers such that
$1\leq\alpha\leq\beta$.
\begin{enumerate}
\item Every homomorphism $f\colon A\to\bigoplus_{j<n}B_j$ is $\lambda$-flat for
some $\lambda\geq 1$.

\item Every $\alpha$-flat homomorphism is $\beta$-flat.

\item Every $\beta$-generic homomorphism is $\alpha$-generic.
\end{enumerate}
\end{lemma}

The definition of canonical generic homomorphisms of
Section~\ref{S:FlGenSp} extends easily to this new context, as
follows.

Let $A$ be a \pss\ space, let $n$ be a natural number, let
$\xf\colon\Idc A\to\two^n$ be a \jzh. So there are unique \jzh s
$\xf_j\colon\Idc A\to\two$, for $j<n$, such that
 \[
 \xf(X)=(\xf_0(X),\ldots,\xf_{n-1}(X)),\qquad\text{for all }
 X\in\Idc A.
 \]
We put $G_j=\Gen(A,\xf_j)$, for all $j<n$, and $G=\bigoplus_{j<n}G_j$.
For any rational number $\lambda\geq1$, we define the \emph{canonical
$\lambda$-generic homomorphism} $f_\lambda\colon A\to G$ by the rule
 \[
 f_\lambda(x)=(f_{0,\lambda}(x),\ldots,f_{n-1,\lambda}(x)),\qquad
 \text{for all }x\in A,
 \]
where $f_{j,\lambda}\colon A\to G_j$ is the
canonical $\lambda$-generic map associated with $A$ and~$\xf_j$, for
all $j<n$. It is trivial, by Lemmas~\ref{L:flflat}
and~\ref{L:flGen} and by Definition~\ref{D:FlGen2}, that
$f_\lambda$ is both $\lambda$-flat and $\lambda$-generic.

Now define a map $\iota\colon\two^n\to\Idc G$ by
 \begin{equation}\label{Eq:Caniota}
\iota(i_0,i_1,\dots,i_{n-1})=\bigoplus_{j\in J}G_j,
 \end{equation}
where we put $J=\setm{j<n}{i_j=1}$. Since the diagram
 \[
{
\def\labelstyle{\displaystyle}
\xymatrix{
\two \ar[r]^(.4){\mathrm{can.}} & \Idc G_j\\
\Idc A\ar[u]^{\xf_j}\ar[ru]_{\Idc f_{j,\lambda}}
}
}
 \]
is commutative for any $j<n$, we get that also the diagram
 \[
{
\def\labelstyle{\displaystyle}
\xymatrix{
\two^n \ar[r]^(.4){\iota} & \Idc G\\
\Idc A\ar[u]^{\xf}\ar[ru]_{\Idc f_\lambda}
}
}
 \]
is commutative. Thus we have proved the following lemma.

\begin{lemma}\label{L:Genlift}
The identity $\iota\circ\xf=\Idc f_\lambda$  holds. Thus the map
$f_\lambda$ lifts $\xf$.
\end{lemma}

The main result of this section is the following analogue of
Lemma~\ref{L:FlGen}(iv).

\begin{lemma}\label{L:FlGen3}
Let $A$, $B$, $G$ be \pss\ spaces,
let $f\colon A\to G$ and $h\colon A\to B$ be positive homomorphisms, and let
$\xg\colon\Idc G\to\Idc B$ be a \jzh\ such that
$\Idc h=\xg\circ\Idc f$. Let $m$, $p$ be the natural numbers
such that $\Idc A\cong\two^m$ and $\Idc G\cong\two^p$, and put
$q=\min\set{2^{m-1},p}$ if $m$, $p>0$, and $q=1$ otherwise.
Let $\lambda\geq1$ be a rational number. If $f$ is
$q\lambda$-generic and $h$ is $\lambda$-flat, then there exists a
homomorphism $g\colon G\to B$ such that $\xg=\Idc g$ and
$h=g\circ f$.
\end{lemma}

\begin{proof}
First of all, by decomposing $B$ as the direct sum of its simple
ideals (see Lemma~\ref{L:DecPss}), it is easy to reduce the
problem to the case where $B$ is \emph{simple}. Furthermore, the problem is
trivial if either $A$ or $G$ is trivial, so suppose from now on that
both $m$ and $p$ are nonzero. Let $A=\bigoplus_{i<m}A_i$ and
$G=\bigoplus_{j<p}G_j$ be the decompositions of $A$ and~$G$ as the direct sums
of their simple ideals. We write
 \[
 f(x)=(f_0(x),\ldots,f_{p-1}(x)),\qquad\text{for all }x\in A,
 \]
with homomorphisms $f_j\colon A\to G_j$, for $j<p$. Moreover, by using
the natural identification of $\Idc G$ with
$\prod_{j<p}\Idc G_j$, we can decompose $\xg$ as
 \[
 \xg(Y_0,\ldots,Y_{p-1})=\bigvee_{j<p}\xg_j(Y_j),\qquad\text{for all }
 (Y_0,\ldots,Y_{p-1})\in\prod_{j<p}\Idc G_j,
 \]
with \jzh s $\xg_j\colon\Idc G_j\to\Idc B$, for $j<p$. Hence
$\Idc h=\bigvee_{j<p}\xh_j$, where we put
$\xh_j=\xg_j\circ\Idc f_j$, for all $j<p$. So $\xh_j$ is a \jzh\
from $\Idc A$ to $\Idc B$.

We define subsets $H$, $H_j$ (for $j<p$) of $m$ as follows:
 \[
 H=\setm{i<m}{h[A_i]\neq\set{0}},\qquad
 H_j=\setm{i<m}{\xh_j(A_i)\neq\set{0}}.
 \]
Put $\CH=\setm{H_j}{j<p}$.
Hence $H=\bigcup_{j<p}H_j=\bigcup\CH$. Furthermore, put
$n_i=|\setm{K\in\CH}{i\in K}|$, for all $i<m$. We make the following
obvious observation:
 \begin{equation}\label{Eq:Ineqni}
 n_i=0,\text{ for all }i\notin H,\qquad\text{and}\qquad
 1\leq n_i\leq q,\quad\text{for all }i\in H.
 \end{equation}
For any $K\in\CH$, we define a homomorphism $h_{(K)}\colon A\to B$ by
its restriction to each~$A_i$, for $i<m$. First, we define the
restriction of $h_{(K)}$ to $A_i$ as equal to zero if $i\notin K$. If
$i\in K$, we put
 \begin{equation}\label{Eq:DefhK}
 h_{(K)}(x)=\dfrac{1}{n_i}h(x),\qquad\text{for all }x\in A_i.
 \end{equation}
Hence, $h_{(K)}[A_i]\neq\set{0}$ if{f} $i\in K$, hence
 \begin{equation}\label{Eq:nabhj}
 \Idc h_{(H_j)}=\xh_j,\qquad\text{for all }j<p.
 \end{equation}
Furthermore, for all $i\in H$ and all $x\in A_i$,
 \[
 \sum_{K\in\CH}h_{(K)}(x)=
 \sum_{\substack{K\in\CH\\ i\in K}}\dfrac{1}{n_i}h(x)=h(x),
 \]
whence we obtain the equality
 \begin{equation}\label{Eq:sumhK}
 \sum_{K\in\CH}h_{(K)}=h.
 \end{equation}
Next, we put $p_j=|\setm{j'<p}{H_j=H_{j'}}|$ and
$h_j=\dfrac{1}{p_j}h_{(H_j)}$, for all $j<p$. Hence $p_j$ is a positive
integer and, by \eqref{Eq:nabhj}, the following holds:
 \begin{equation}\label{Eq:nabhj2}
 \Idc h_j=\xh_j,\qquad\text{for all }j<p.
 \end{equation}
Moreover, $\sum_{j<p}h_j=\sum_{K\in\CH}h_{(K)}=h$ (the last equality follows from
\eqref{Eq:sumhK}). So we have obtained
 \begin{equation}\label{Eq:sumhj}
 \sum_{j<p}h_j=h.
 \end{equation}

Now we prove the following crucial claim:

\begin{claim}\label{Cl:lql}
The homomorphism $h_j$ is $q\lambda$-flat, for all $j<p$.
\end{claim}

\begin{cproof}
By the definition of $h_j$, it suffices to prove that $h_{(K)}$ is
$q\lambda$-flat, for all $K\in\CH$. For all $i\in K$,
$h_{(K)}(1_{A_i})=\dfrac{1}{n_i}h(1_{A_i})$, so, for all $i$, $j\in K$,
 \begin{align*}
 h_{(K)}(1_{A_j})=\dfrac{1}{n_j}h(1_{A_j})&\lear
 \dfrac{\lambda}{n_j}h(1_{A_i})&&
 (\text{because }h\text{ is }\lambda\text{-flat})\\
 &=\dfrac{n_i\lambda}{n_j}h_{(K)}(1_{A_i})\\
 &\leq n_i\lambda h_{(K)}(1_{A_i})\\
 &\leq q\lambda h_{(K)}(1_{A_i})&&(\text{by \eqref{Eq:Ineqni}}),
 \end{align*}
which concludes the proof of the claim.
\end{cproof}

Since $f$ is $q\lambda$-generic, it follows from Lemma~\ref{L:FlGen}
that for all $j<p$, there exists a homomorphism $g_j\colon G_j\to B$
such that $\Idc g_j=\xg_j$ and $g_j\circ f_j=h_j$. We define
$g\colon G\to B$ by the rule
 \[
 g(y_0,\ldots,y_{p-1})=\sum_{j<p}g_j(y_j),\qquad\text{for all }
 (y_0,\ldots,y_{p-1})\in B.
 \]
Then $\Idc g=\xg$ and, by \eqref{Eq:sumhj}, $g\circ f=h$.
\end{proof}

\begin{example}
The following example shows that one cannot replace $q$ by
$1$ in the statement of Lemma~\ref{L:FlGen3}. Let $h\colon\QQ^2\to\QQ$
be defined by the rule $h(x,y)=x+y$ (for all $x$, $y\in\QQ$), and let
$f\colon\QQ^2\to A\oplus B$ be any $1$-generic lifting of the map
$\xf\colon\two^2\to\two^2$ defined by the rule $\xf(\bx,\by)=(\bx\vee\by,\by)$,
with $A$ and $B$ being simple \pss\ spaces. Let $\xg\colon\two^2\to\two$ be defined
by the rule $\xg(\bx,\by)=\bx\vee\by$. By identifying $\Idc(A\oplus B)$ with
$\two^2$ and $\Idc\QQ$ with $\two$, the map $\xg$ defines a \jzh\
from $\Idc(A\oplus B)$ to $\Idc\QQ$, and
$\Idc h=\xg\circ\Idc f$. Obviously, $h$ is $1$-flat. However, we
prove that there exists no lifting $g$ of $\xg$ such that $h=g\circ f$. The
situation is summarized on the diagrams below:
 \[
{
\def\labelstyle{\displaystyle}
\xymatrix{
 & \Idc\QQ & & & \QQ\\
 \Idc(A\oplus B)\ar[ur]^{\xg} & & & A\oplus B\ar@{-->}[ur]^g\\
 & \Idc(\QQ\oplus\QQ)\ar[lu]^{\Idc f}\ar[-2,0]_{\Idc h} & & &
 \QQ\oplus\QQ\ar[lu]^f\ar[-2,0]_h
}
}
 \]
There are $a$, $c\in A^{++}$ and $b\in B^{++}$ such that
 \[
 f(x,y)=(xa+yc,yb),\qquad\text{for all }x,\,y\in\QQ.
 \]
{}From the fact that $f$ is $1$-flat follows that $a\asar_1c$ in $A$.
Furthermore, there are positive, nonzero homomorphisms
$u\colon A\to\QQ$ and $v\colon B\to\QQ$ such that
 \[
 g(x,y)=u(x)+v(y),\qquad\text{for all }x\in A\text{ and }y\in B.
 \]
{}From $a\asar_1c$ follows that $u(a)\asar_1u(c)$, thus, since $u(a)$,
$u(c)\in\QQ^+$, $u(a)=u(c)$. However, by applying $h=g\circ f$ to $(1,0)$ and
$(0,1)$ respectively, we obtain that $u(a)=1$ and $u(c)+v(b)=1$, hence
$0<u(c)<u(a)=u(c)$, a contradiction.
\end{example}

We end this section with an easy lifting result:

\begin{lemma}\label{L:RevLift}
Let $E$ be a \pss\ space, let $m$ be a natural number, let
$\xf\colon\two^m\to\Idc E$ be a \jzh. Then there exists a positive homomorphism
$f\colon\QQ^m\to E$ such that, if $\iota\colon\two^m\to\Idc\QQ^m$ denotes the
canonical isomorphism, the equality $\xf=\Idc f\circ\iota$ holds.
\end{lemma}

\begin{proof}
Let $E=\bigoplus_{j<n}E_j$ be the decomposition of $E$ as a direct sum of simple
\pss\ spaces. Denote by $(\xe_i)_{i<m}$ (resp., $(e_i)_{i<m}$) the canonical
basis of $\two^m$ (resp., $\QQ^m$). For all $i<m$, there exists a subset $J_i$ of
$\set{0,1,\ldots,n-1}$ such that $\xf(\xe_i)=\bigoplus_{j\in J_i}E_j$. Let
$f\colon\QQ^m\to E$ be the linear map defined by $f(e_i)=\sum_{j\in J_i}1_{E_j}$,
for all $i<m$. Then, for all $i<m$,
 \[
 (\Idc f\circ\iota)(\xe_i)=(\Idc f)(\QQ e_i)=E(f(e_i))=\bigoplus_{j\in J_i}E_i
 =\xf(\xe_i),
 \]
whence $\Idc f\circ\iota=\xf$.
\end{proof}

\section{\Pss\ liftings of dismantlable diagrams}\label{S:DismLift}

The existence of $\lambda$-generic mappings enables us to construct liftings
of a large class of finite diagrams with respect to the functor $\Idc$. The
following definition is due to K.\,A. Baker, P.\,C. Fishburn, and F.\,S. Roberts
\cite{BFR}.

\begin{definition}\label{D:Dismantl}
A finite lattice $L$ of cardinality $n$ is called \emph{dismantlable} if there
is a chain $L_1\subset L_2\subset\cdots\subset L_n=L$ of sublattices of $L$
such that $|L_i|=i$ for every $i=1$, $2$, \dots, $n$.
\end{definition}

We shall use an extension of Definition~\ref{D:Dismantl} to partially ordered
sets. \emph{We shall consider $\es$ as a partially ordered set}.
For a finite partially ordered set $P$, an element $x$ of $P$ is
\emph{doubly-irreducible}, if $x$ has at most one upper cover and at most one
lower cover. We denote by $\Irr P$ the set of all doubly-irreducible elements of
$P$. For a subset $Q$ of $P$, let $Q\lessdot P$ be the statement that
$Q=P\setminus\set{x}$ for some $x\in\Irr P$.

\begin{definition}\label{D:Dismantpo}
A finite partially ordered set $P$ of cardinality $n$ is \emph{dismantlable}, if
there exists a chain $\es=P_0\lessdot P_1\lessdot\cdots\lessdot P_n=P$ of subsets
of $P$.
\end{definition}

Of course, a finite lattice is dismantlable if{f} it is dismantlable as a
partially ordered set. On the other hand, every finite bounded
dismantlable partially ordered set is a lattice as observed
in D. Kelly and I. Rival \cite{KR2}.
Furthermore, a nonempty partially ordered set $P$ is
dismantlable if{f} there exists $x\in\Irr P$ such that $P\setminus\set{x}$ is
dismantlable.

There is a large supply of finite dismantlable lattices. For example,
the following result was proved in \cite{BFR}; another proof can be found in
D. Kelly and I. Rival \cite{KR2}.

\begin{theorem}\label{T:planar}
Every finite planar lattice is dismantlable.
\end{theorem}

There are also non-planar modular dismantlable lattices, like the one
diagrammed on Figure~1. On the other hand, every distributive dismantlable
lattice is planar, see D. Kelly and I. Rival \cite{KR1}.

\begin{figure}[ht]
\includegraphics{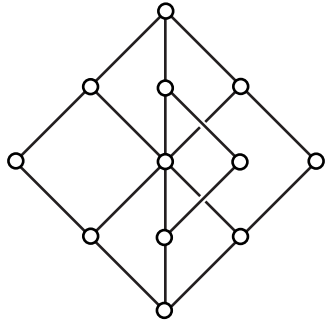}
\caption{A dismantlable, modular, non-planar lattice}
\end{figure}

The main result of this section is the following.

\begin{theorem}\label{T:dislift}
Every diagram $\Phi\colon P\to\SSfb$ of finite Boolean
semilattices indexed by a finite dismantlable partially ordered set $P$ has a
\pss\ lifting $\Psi\colon P\to\EEps$.
\end{theorem}

\begin{proof}
We shall proceed by induction on the cardinality of $P$.
We may also assume without loss of generality that for every $x\in P$ we have
$\Phi(x)=\two^q$ for some $q\in\ZZ^+$. If $P=\es$ then the result is trivial.

Now assume that $P$ is nonempty. So there exists $x\in\Irr P$ such that
$Q=P\setminus\set{x}$ is dismantlable. By the induction hypothesis, the
restriction $\Phi'$ of $\Phi$ to $Q$ has a \pss\ lifting
$\Psi'\colon Q\to\EEps$. So there are isomorphisms
$\iota_t\colon\Phi(t)\to\Idc(\Psi'(t))$, for
$t\in Q$, such that the following holds:
 \begin{equation}\label{Eq:CommQ}
 \iota_t\circ\Phi(\eps_{s,t})=\Idc(\Psi'(\eps_{s,t}))\circ\iota_s,
 \qquad\text{for all }s\leq t\text{ in }Q.
 \end{equation}
Let $n\in\ZZ^+$ such that $\Phi(x)=\two^n$.

Suppose first that $x$ is comparable to no element of $Q$. Then, in order to lift
$\Phi$, it is sufficient to adjoin one lifting of $\Phi(x)$ to the diagram
$\Psi'$; take $\Psi(x)=\QQ^n$.

So, from now on, suppose that $x$ can be compared to at least one element of $Q$.
We denote by $u$ (resp., $v$) the unique lower cover (resp., upper cover) of $x$
in $P$ if it exists. By assumption, either $u$ or $v$ exists.

We start with the case where $x$ is a minimal element of $P$. Then $v$ exists,
and, by Lemma~\ref{L:RevLift}, there exists a positive homomorphism
$f\colon\QQ^n\to\Psi'(v)$ such that, if $\iota_x\colon\two^n\to\Idc\QQ^n$ denotes
the canonical isomorphism, the equality
 \begin{equation}\label{Eq:easyiot}
 \iota_v\circ\Phi(\eps_{x,v})=\Idc f\circ\iota_x
 \end{equation}
holds. We define $\Psi(x)=\QQ^n$, and, for every $t\in Q$ such that $x<t$, that
is, $v\leq t$, we put $\Psi(\eps_{x,t})=\Psi'(\eps_{v,t})\circ f$. It follows
from \eqref{Eq:CommQ} that the right half of the following diagram commutes,
 \[
{
\def\labelstyle{\displaystyle}
\xymatrix{
 \Idc(\Psi(x))=\Idc(\QQ^n)\ar[rr]^{\Idc f} & &
 \Idc(\Psi'(v))\ar[rr]^{\Idc(\Psi'(\eps_{v,t}))} & & \Idc(\Psi'(t))\\
 \Phi(x)=\two^n\ar[u]^{\iota_x}\ar[rr]_{\Phi(\eps_{x,v})}
 & & \Phi(v)\ar[u]_{\iota_v}\ar[rr]_{\Phi(\eps_{v,t})} & &
 \Phi(t)\ar[u]_{\iota_t}
}
}
 \]
while the commutation of the left half of this diagram follows from
\eqref{Eq:easyiot}. This settles the case where $x$ is a minimal element of $P$.

So, suppose, from now on, that $x$ is not minimal in $P$, so that $u$ exists. Let
$m\in\ZZ^+$ such that $\Phi(u)=\two^m$. Observe that
$\xf=\Phi(\eps_{u,x})\circ\iota_u^{-1}$ is a \jzh\ from $\Idc(\Psi'(u))$ to
$\two^n$. Put $G=\Gen(\Psi'(u),\xf)$ (see Section~\ref{S:FlGen}), let
$\iota_x\colon\two^n\to\Idc G$ be the isomorphism given by \eqref{Eq:Caniota}.
For any rational number $\mu\geq1$, denote by $f_\mu\colon\Psi'(u)\to G$ the
canonical $\mu$-generic lifting of $\xf$, see
Section~\ref{S:FlGen}. Hence, by Lemma~\ref{L:Genlift}, the following diagram
commutes,
 \[
{
\def\labelstyle{\displaystyle}
\xymatrix{
\two^n \ar[r]^(.4){\iota_x} & \Idc G\\
\Idc(\Psi'(u))\ar[u]^{\xf}\ar[ru]_{\Idc f_\mu}
}
}
 \]
which means that the following equality holds:
 \begin{equation}\label{Eq:xfmu}
 \Idc f_\mu\circ\iota_u=\iota_x\circ\Phi(\eps_{u,x}).
 \end{equation}
In case $x$ is maximal in $P$, set $\mu=1$. In case $x$ is not maximal in $P$,
there exists, by Lemma~\ref{L:FlGen2}(i) a rational number $\lambda\geq1$ such
that $\Psi'(\eps_{u,v})$ is $\lambda$-flat.
Put $q=\min\set{2^{m-1},n}$ if $m$, $n>0$, $q=1$ otherwise, and put
$\mu=q\lambda$. In both cases, we define $\Psi(x)=G$ and
$\Psi(\eps_{u,x})=f_\mu$, the canonical
$\mu$-generic lifting of $\xf$. Hence for every $t\leq u$ we have to set
$\Psi(\eps_{t,x})=f_\mu\circ\Psi'(\eps_{t,u})$, and we need to check that the
following diagram commutes:
 \[
{
\def\labelstyle{\displaystyle}
\xymatrix{
 \Idc(\Psi'(t))\ar[rr]^{\Idc(\Psi'(\eps_{t,u}))} & &
 \Idc(\Psi'(u))\ar[rr]^{\Idc f_\mu} & & \Idc(\Psi(x))=\Idc G\\
 \Phi(t)\ar[u]^{\iota_t}\ar[rr]_{\Phi(\eps_{t,u})}
 & & \Phi(u)\ar[u]_{\iota_u}\ar[rr]_{\Phi(\eps_{u,x})} & &
 \Phi(x)\ar[u]_{\iota_x}
}
}
 \]
The commutativity of the left half of this diagram follows from \eqref{Eq:CommQ},
while the commutativity of the right half follows from \eqref{Eq:xfmu}.
In particular, this completes the proof of the induction step if $x$ is the largest
element of $P$.

If $x$ has a (unique) upper cover $v$ we continue as follows. We first observe
that $\xg=\iota_v\circ\Phi(\eps_{x,v})\circ\iota_x^{-1}$ is a \jzh\ from $\Idc G$
to $\Idc(\Psi'(v))$.  Since $\Psi'(\eps_{u,v})$ is
$\lambda$-flat and $\Psi'(\eps_{u,x})=f_\mu$ is $q\lambda$-generic,
Lemma~\ref{L:FlGen3} can be applied to the two following commutative diagrams:
 \[
{
\def\labelstyle{\displaystyle}
\xymatrix{
 \Psi'(v) &  &  &  & \Idc\Psi'(v) & \\
 & \Psi(x)=G &  &  &  &
 \Idc G\ar[ul]_{\xg}\\
 \Psi'(u)\ar[uu]^{\Psi'(\eps_{u,v})}\ar[ur]_{f_\mu} &  &  &  &
 \Idc(\Psi'(u))\ar[uu]^{\Idc(\Psi'(\eps_{u,v}))}
 \ar[ur]_{\iota_x\circ\xf}
}
}
 \]
Thus there exists a positive homomorphism $g\colon G\to\Psi'(v)$ such that
$\Idc g=\xg$ and $g\circ f_\mu=\Psi'(\eps_{u,v})$. We define $\Psi(\eps_{x,v})=g$,
thus for every $t\geq v$ we have to set
$\Psi(\eps_{x,t})=\Psi'(\eps_{v,t})\circ g$, and we need to check that the
following diagram commutes:
 \[
{
\def\labelstyle{\displaystyle}
\xymatrix{
 \Idc(\Psi'(u))\ar[rr]^{\Idc f_\mu} & &
 \Idc G\ar[rr]^{\xg=\Idc g} & & \Idc(\Psi'(v))\ar[rr]^{\Idc(\Psi'(\eps_{v,t}))}
 & & \Idc(\Psi'(t))\\
 \Phi(u)\ar[u]^{\iota_u}\ar[rr]_{\Phi(\eps_{u,x})}
 & & \Phi(x)\ar[u]_{\iota_x}\ar[rr]_{\Phi(\eps_{x,v})} & &
 \Phi(v)\ar[u]_{\iota_v}\ar[rr]_{\Phi(\eps_{v,t})} & & \Phi(t)\ar[u]_{\iota_t}
}
}
 \]
The commutativity of the right third of this diagram follows from
\eqref{Eq:CommQ}, the commutativity of the left third follows from
\eqref{Eq:xfmu}, and the commutativity of the middle third of this diagram
follows from the fact that $\Idc g=\xg$ and the definition of~$\xg$. This
completes the verification that the extension $\Psi$ of $\Psi'$ thus defined is
as required.
\end{proof}

\section{Lifting \jzh s between countable \jzs s}\label{S:CtbleSem}

The main result of this section is a wide generalization of
Lemma~\ref{L:RevLift}:

\begin{theorem}\label{T:Nab1-1}
Let $S$ be a countable distributive \jzs, let $H$ be a countable dimension
vector space, let $\xf\colon S\to\Idc H$ be a \jzh. Then $\xf$ can be lifted,
that is, there are a countable dimension vector space $G$, a positive
homomorphism $f\colon G\to H$, and an isomorphism $\alpha\colon S\to\Idc G$ such
that $\xf=(\Idc f)\circ\alpha$.
\end{theorem}

\begin{proof}
It follows from E.\,G. Effros, D.\,E. Handelman, and C.-L. Shen~\cite{EHS80} (see
also \cite[Theorem~3.19]{Gpoag}) that $H$ is the direct limit of a sequence of
simplicial (partially ordered) groups, that is, groups of the form $\ZZ^n$,
for $n\in\ZZ^+$. Since $H$ is a vector space, the isomorphism
$H\cong H\otimes\QQ$ implies that $H$ is, in fact, the direct limit of a
sequence of simplicial \emph{vector spaces}: more specifically,
$H=\varinjlim_{i\in\ZZ^+}H_i$, where $H_i=\QQ^{n_i}$ for some $n_i\in\ZZ^+$,
with transition maps $t_i\colon H_i\to H_{i+1}$ and limiting maps
$t'_i\colon H_i\to H$, both being positive homomorphisms, for all $i\in\ZZ^+$.

Furthermore, it follows from \cite[Corollary~6.7(i)]{GoWe} that $S$ is the
direct limit of a countable sequence of finite Boolean \jzs s, say,
$S=\varinjlim_{i\in\ZZ^+}\two^{m_i}$, with $m_i\in\ZZ^+$, transition maps
$\xs_i\colon\two^{m_i}\to\two^{m_{i+1}}$ and limiting maps
$\xs'_i\colon\two^{m_i}\to S$, both being \jzh s, for all $i\in\ZZ^+$. The
situation may be visualized by the following commutative diagrams:
 \[
{
\def\labelstyle{\displaystyle}
\xymatrix{
 \two^{m_{i+1}}\ar[dr]^{\xs'_{i+1}} &  &  &  &
 H_{i+1}\ar[dr]^{t'_{i+1}} & \\
 \two^{m_i}\ar[u]^{\xs_i}\ar[r]_{\xs'_i} & S &  &  &
 H_i\ar[u]^{t_i}\ar[r]_{t'_i} & H
}
}
 \]
For all $i\in\ZZ^+$, the \jzh\
$\xf\circ\xs'_i\colon\two^{m_i}\to\Idc H\cong\varinjlim_{j\in\ZZ^+}\Idc H_j$
factors through $\Idc H_j$ for some $j\in\ZZ^+$. Hence, by possibly replacing the
sequence $(H_j)_{j\in\ZZ^+}$ by some subsequence, we can suppose that $j=i$,
allowing to factor $\xf$ through a \jzh\ $\xf_i\colon\two^{m_i}\to\Idc H_i$. The
situation may be visualized by the following commutative diagram, where the two
parallel horizontal sequences are direct limits:
 \begin{equation}\label{Eq:IdcHitwomi}
{
\def\labelstyle{\displaystyle}
\xymatrix{
 \Idc H_0\ar[r]^{\Idc t_0}\ar@/^4pc/[rrrr]|-{\Idc t'_0}
 & \Idc H_1\ar[r]^(.7){\Idc t_1}\ar@/^2pc/[rrr]|-{\Idc t'_1}
 & {}\ar@{--}[r]
 & {}\ar[r] & \Idc H\\
 \two^{m_0}\ar[u]^{\xf_0}\ar[r]_{\xs_0}\ar@/_4pc/[rrrr]|-{\xs'_0}
 & \two^{m_1}\ar[r]_(.6){\xs_1}\ar[u]^{\xf_1}\ar@/_2pc/[rrr]|-{\xs'_1}
 & {}\ar@{--}[r]
 & {}\ar[r] & S\ar[u]_{\xf}
}
}
 \end{equation}
Now we construct, for $i\in\ZZ^+$, a certain \pss\ space $G_i$, together with an
isomorphism $\alpha_i\colon\two^{m_i}\to\Idc G_i$ and a positive homomorphism
$f_i\colon G_i\to H_i$ such that the following equality holds:
 \begin{equation}\label{Eq:fiIdcfialph}
 \xf_i=(\Idc f_i)\circ\alpha_i.
 \end{equation}
For $i=0$, we pick $G_0$, $f_0$, and $\alpha_0$ satisfying
\eqref{Eq:fiIdcfialph}; their existence is ensured by Lemma~\ref{L:RevLift}. Now
suppose $G_i$, $f_i$, and $\alpha_i$ constructed satisfying
\eqref{Eq:fiIdcfialph} above, we construct $G_{i+1}$, $f_{i+1}$, and
$\alpha_{i+1}$. First, by Lemma~\ref{L:FlGen2}(i), there exists a rational
number $\lambda\geq1$ such that the homomorphism
$t_i\circ f_i\colon G_i\to H_{i+1}$ is $\lambda$-flat. We put
$q=\min\set{2^{m_i-1},n_{i+1}}$ if $m_i$, $n_{i+1}>0$
and $q=1$ otherwise. Furthermore, we put
$G_{i+1}=\Gen(G_i,\xs_i\circ\alpha_i^{-1})$, and we let
$s_i\colon G_i\to G_{i+1}$ be the canonical $q\lambda$-generic lifting of
$\xs_i\circ\alpha_i^{-1}$. In particular, Lemma~\ref{L:Genlift} yields us an
isomorphism $\alpha_{i+1}\colon\two^{m_{i+1}}\to\Idc G_{i+1}$ such that the
following equality holds:
 \begin{equation}\label{Eq:DefAlphi+1}
 \Idc s_i=\alpha_{i+1}\circ\xs_i\circ\alpha_i^{-1}.
 \end{equation}
It follows from Lemma~\ref{L:FlGen3} that there exists a positive
homomorphism
$f_{i+1}\colon G_{i+1}\to H_{i+1}$ such that
$\Idc f_{i+1}=\xf_{i+1}\circ\alpha_{i+1}^{-1}$ and
$f_{i+1}\circ s_i=t_i\circ f_i$. The situation may be visualized by the
following commutative diagrams:
 \[
{
\def\labelstyle{\displaystyle}
\xymatrix{
 H_{i+1} &  &  &  & \Idc H_{i+1} & \\
 & G_{i+1}\ar[ul]_{f_{i+1}} &  &  &  &
 \Idc G_{i+1}\ar[ul]_{\xf_{i+1}\circ\alpha_{i+1}^{-1}}\\
 G_i\ar[uu]^{t_i\circ f_i}\ar[ur]_{s_i} &  &  &  &
 \Idc G_i\ar[uu]^{\Idc(t_i\circ f_i)}
 \ar[ur]_{\alpha_{i+1}\circ\xs_i\circ\alpha_i^{-1}}
}
}
 \]
Hence the following diagrams are commutative:
 \begin{equation}\label{Eq:fifi+1alph}
{
\def\labelstyle{\displaystyle}
\xymatrix{
 H_i\ar[r]^{t_i} & H_{i+1} &  \Idc G_i\ar[r]^{\Idc s_i} & \Idc G_{i+1}
 & \Idc G_{i+1}\ar[dr]^{\Idc f_{i+1}} & \\
 G_i\ar[u]^{f_i}\ar[r]_{s_i} & G_{i+1}\ar[u]_{f_{i+1}} &
 \two^{m_i}\ar[u]^{\alpha_i}\ar[r]_{\xs_i} & \two^{m_{i+1}}\ar[u]_{\alpha_{i+1}}
 & \two^{m_{i+1}}\ar[u]^{\alpha_{i+1}}\ar[r]_{\xf_{i+1}} & \Idc H_{i+1}
}
}
 \end{equation}
We let $G$ be the direct limit of the sequence $(G_i)_{i\in\ZZ^+}$, with
transition maps
$s_i\colon G_i\to G_{i+1}$ and limiting maps
$s'_i\colon G_i\to G$, for all $i\in\ZZ^+$. Hence $G$ is a dimension vector
space. Furthermore, since the functor $\Idc$ preserves direct limits, the
commutativity of the left diagram and the middle diagram in
\eqref{Eq:fifi+1alph} imply the existence
of a positive homomorphism $f\colon G\to H$ and an isomorphism
$\alpha\colon S\to\Idc G$ such that the following infinite diagrams are
commutative:
 \begin{equation}\label{Eq:G,Hrel}
{
\def\labelstyle{\displaystyle}
\xymatrix{
 H_0\ar[r]^{t_0}
 & H_1\ar[r]^(.6){t_1}
 & {}\ar@{--}[r]
 & {}\ar[r] & H\\
 G_0\ar[u]^{f_0}\ar[r]_{s_0}
 & G_1\ar[u]^{f_1}\ar[r]_(.6){s_1}
 & {}\ar@{--}[r]
 & {}\ar[r] & G\ar[u]_{f}
}
}
\end{equation}
and
 \begin{equation}\label{Eq:IdcG,Srel}
{
\def\labelstyle{\displaystyle}
\xymatrix{
 \Idc G_0\ar[r]^{\Idc s_0}
 & \Idc G_1\ar[r]^(.7){\Idc s_1}
 & {}\ar@{--}[r]
 & {}\ar[r] & \Idc G\\
 \two^{m_0}\ar[u]^{\alpha_0}\ar[r]_{\xs_0}
 & \two^{m_1}\ar[u]^{\alpha_1}\ar[r]_(.6){\xs_1}
 & {}\ar@{--}[r]
 & {}\ar[r] & S\ar[u]_{\alpha}
}
}
\end{equation}
(For sake of clarity, we do not represent on these diagrams the transition maps
corresponding to the direct limits in the rows of \eqref{Eq:G,Hrel},
\eqref{Eq:IdcG,Srel}.) To verify the equality $(\Idc f)\circ\alpha=\xf$, it
suffices to verify that
$(\Idc f_i)\circ\alpha_i=\xf_i$ holds, for all $i\in\ZZ^+$, which is
exactly~\eqref{Eq:fiIdcfialph}.
\end{proof}

In particular, it follows from Theorem~\ref{T:Nab1-1} that every countable
distributive \jzs\ is isomorphic to $\Idc G$ for some dimension vector space
$G$, but this is of course much easier to prove directly, see
\cite[Theorem~5.2]{GoWe}.
As an immediate consequence of this and Theorem~\ref{T:Nab1-1}, we record the
following:

\begin{corollary}\label{C:Nab1-1}
Let $S$ and $T$ be countable distributive \jzs s, let
$\xf\colon S\to T$ be a \jzh. Then $\xf$ can be lifted, that is, there are
countable dimension vector spaces $G$ and $H$ together with a positive
homomorphism $f\colon G\to H$ and isomorphisms
$\alpha\colon S\to\Idc G$ and $\beta\colon T\to\Idc H$ such that
$(\Idc f)\circ\alpha=\beta\circ\xf$.
\end{corollary}

Corollary~\ref{C:Nab1-1} cannot be extended to uncountable semilattices, see
Example~\ref{Ex:al1diagnolift}.

In order to obtain a ring-theoretical version of Corollary~\ref{C:Nab1-1}, we
shall need the following more precise well-known form of the isomorphism
\eqref{Eq:IdG=IdR}, which is essentially contained
in~\cite[Theorem~15.20]{Gvnrr}:

\begin{lemma}\label{L:IdG=IdR}
Let $K$ be a field. For a unital locally matricial $K$-algebra $R$, we let
$\eta_R\colon\Idc(K_0(R))\to\Idc R$ be the map defined by the rule
 \[
 \eta_R(I)=\setm{x\in R}{[xR]\in I},\quad\text{for all }I\in\Idc(K_0(R))
 \]
\pup{we let $[xR]$ denote the isomorphism type of $xR$}. Then $\eta_R$ is a
semilattice isomorphism and the rule $R\mapsto\eta_R$ defines a natural
equivalence from the functor $\Idc\circ K_0$ \pup{$\Idc$ is meant on dimension
groups} to the functor
$\Idc$ \pup{defined on locally matricial $K$-algebras}. This means that for any
unital locally matricial $K$-algebras $A$ and $B$ and any unital $K$-algebra
homomorphism $f\colon A\to B$, the following diagram is commutative:
 \[
{
\def\labelstyle{\displaystyle}
\xymatrix{
 \Idc A\ar[rr]^{\Idc f} & & \Idc B\\
 \Idc(K_0(A))\ar[rr]_{\Idc(K_0(f))}\ar[u]^{\eta_A}
 & & \Idc(K_0(B))\ar[u]_{\eta_B}
}
}
 \]
\end{lemma}

We can now state the ring-theoretical analogue of Corollary~\ref{C:Nab1-1}:

\begin{corollary}\label{C:Bab1-1ring}
Let $S$ and $T$ be countable distributive \jzs s \pup{resp., countable
distributive \jzs s with unit}, let $\xf\colon S\to T$ be a \jzh\
\pup{resp., a \jzuh}, let $K$ be a field. Then there are
countably dimensional locally matricial \pup{resp., countably dimensional locally
matricial unital} $K$-algebras $A$ and $B$, a $K$-algebra homomorphism
\pup{resp., a unital $K$-algebra homomorphism} $f\colon A\to B$, and isomorphisms
$\alpha\colon S\to\Idc A$ and $\beta\colon T\to\Idc B$ such that
$(\Idc f)\circ\alpha=\beta\circ\xf$.
\end{corollary}

\begin{proof}
We first deal with the case of unital homomorphisms, so both $S$ and
$T$ have unit and $\xf(1)=1$. By Corollary~\ref{C:Nab1-1}, there are
countable dimension vector spaces $G$ and $H$ together with a positive
homomorphism $f'\colon G\to H$ and isomorphisms
$\alpha'\colon S\to\Idc G$ and $\beta'\colon T\to\Idc H$ such that
$(\Idc f')\circ\alpha'=\beta'\circ\xf$. Since both $S$ and $T$ have a largest
element, both $G$ and $H$ have an order-unit, that we will denote respectively
by $1_G$ and~$1_H$. Moreover, from $\xf(1)=1$ follows that we may take
$f'(1_G)=1_H$. By the Elliott, Effros, Handelman, and Shen Theorem (see
\cite{Elli76,EHS80,Gvnrr}), there are countably dimensional, locally matricial,
unital $K$-algebras $A$,~$B$ and (normalized) isomorphisms
$\alpha''\colon(G,1_G)\to(K_0(A),[A])$ and
$\beta''\colon(H,1_H)\to(K_0(B),[B])$. By \cite[Lemma~1.3]{GoHa86}, there exists
a unital $K$-algebra homomorphism $f\colon A\to B$ such that
$K_0(f)=\beta''\circ f'\circ{\alpha''}^{-1}$. Hence the following diagram
commutes:
 \[
{
\def\labelstyle{\displaystyle}
\xymatrix{
 \Idc(K_0(A))\ar[rr]^{\Idc(K_0(f))} & & \Idc(K_0(B))\\
 S\ar[rr]_{\xf}\ar[u]^{(\Idc\alpha'')\circ\alpha'}
 & & T\ar[u]_{(\Idc\beta'')\circ\beta'}
}
}
 \]
By Lemma~\ref{L:IdG=IdR}, the following diagram commutes,
 \[
{
\def\labelstyle{\displaystyle}
\xymatrix{
 \Idc A\ar[r]^{\Idc f} & \Idc B\\
 S\ar[r]_{\xf}\ar[u]^{\alpha}
 & T\ar[u]_{\beta}
}
}
 \]
where we put $\alpha=\eta_A\circ(\Idc\alpha'')\circ\alpha'$ and
$\beta=\eta_B\circ(\Idc\beta'')\circ\beta'$. This solves the case with unit.

In the case without unit, we put $S'=S\cup\set{1}$, $T'=T\cup\set{1}$ (where $1$
is the largest element), and we extend $\xf$ to a map $\xf'\colon S'\to T'$ by
putting
$\xf'(1)=1$. We let $f'\colon A'\to B'$, together with isomorphisms
$\alpha'\colon S'\to\Idc A'$ and $\beta'\colon T'\to\Idc B'$, lift $\xf'$, by
using the result of the case with unit. Then we put
 \begin{align*}
 A&=\setm{x\in A'}{{\alpha'}^{-1}(A'xA')\in S},\\
 B&=\setm{y\in B'}{{\beta'}^{-1}(B'yB')\in T},
 \end{align*}
so $A$ (resp., $B$) is an ideal of $A'$ (resp., $B'$).
Then we let $f$ be the restriction of $f'$ from $A$ to $B$ and
$\alpha\colon S\to\Idc A$ (resp., $\beta\colon T\to\Idc B$) be the restrictions
of $\alpha'$ and $\beta'$ to $S$ and $T$, respectively.
\end{proof}

Every locally matricial algebra $R$ is a von~Neumann regular ring, and then
$\Idc R$ is isomorphic to $\Conc\LL(R)$, the semilattice of compact congruences
of the sectionally complemented, modular lattice $\LL(R)$ of all principal right
ideals of $R$, see \cite[Corollary~4.4]{Wehr99}. (A lattice $L$ with zero is
\emph{sectionally complemented}, if for all $a\leq b$ in $L$, there exists
$x\in L$ such that $a\wedge x=0$ and $a\vee x=b$.) The proof of
\cite[Corollary~4.4]{Wehr99} does not involve the unit, and, for a finite
field $K$, every locally matricial $K$-algebra $R$ is locally finite, thus so is
the lattice $\LL(R)$. Thus similar methods as those used in the proof of
Corollary~\ref{C:Bab1-1ring} yield easily the following lattice-theoretical
analogue of Corollary~\ref{C:Bab1-1ring}:

\begin{corollary}\label{C:Bab1-1latt}
Let $S$ and $T$ be countable distributive \jzs s, let $\xf\colon S\to T$ be a
\jzh. Then there are locally finite, sectionally complemented, modular lattices
$L$ and $M$ together with a \jzlh\ $f\colon L\to M$ and isomorphisms
$\alpha\colon S\to\Conc L$ and $\beta\colon T\to\Conc M$ such that
$(\Conc f)\circ\alpha=\beta\circ\xf$. Furthermore, if both $S$ and $T$ have a
unit, then one may take both $L$ and $M$ with unit and $f(1)=1$.
\end{corollary}

The lattices $L$ and $M$ are more than just locally finite and sectionally
complemented, for example, every finite subset of $L$ is contained in a finite,
sectionally complemented, modular sublattice of $L$.

\section{$1$-dimensional lifting fails for the $\Idc$ functor and
\jzs\ embeddings}

The following example is a modification of an earlier example of the second
author, and it has been communicated to the authors by P. R\r{u}\v{z}i\v{c}ka.

\begin{examplepf}\label{Ex:nabla1fails}
There exists a dimension group $G$ together with a \jz-embedding
$\mu\colon\Idc G\hookrightarrow\two^2$ that cannot be
lifted. More precisely, there are no {\poag}~$H$, no positive homomorphism
$f\colon G\to H$, and no isomorphism $\alpha\colon\Idc H\to\two^2$ such that
$\mu=\alpha\circ\Idc f$.
\end{examplepf}

\begin{proof}
Let $G=\ZZ\times_{\mathrm{lex}}\ZZ$ be the lexicographical product of the chain
$\ZZ$ of integers with itself, put
$a=(0,1)$ and $b=(1,0)$. Hence $na\leq b$, for all $n\in\ZZ^+$, which we
write $a\ll b$. So $\Idc G\cong\three$, the three-element chain. Denote
further by $x$, $y$ the atoms of $\two^2$.

Let $\mu\colon\Idc G\to\two^2$ be defined by the rule $\mu(G(0))=0$,
$\mu(G(a))=x$, and $\mu(G(b))=1$. Observe that $\mu$ is a
\jz-embedding. Suppose that there are $H$, $f$, and $\alpha$ as in the
statement above.

Put $a'=f(a)$ and $b'=f(b)$.
{}From $a\ll b$ follows that $a'\ll b'$. {}From $\alpha(H(b'))=1=x\vee y$
follows that there are $c$, $d\in H^+$ such that
$\alpha(H(c))=x$,
$\alpha(H(d))=y$, and $b'\leq c+d$. Since $\alpha(H(a'))=x=\alpha(H(c))$ and
$\alpha$ is an isomorphism, the relation $a'\as c$ holds, thus there exists
$n\in\NN$ such that $c\leq na'$. Put $u=b'-na'$. Since $0<a'\ll b'$, we get
$0<a'\ll u$ and $H(u)=H(b')$; thus $\alpha(H(u))=1$.
But from $c\leq na'$ we get $u=b'-na'\leq b'-c\leq d$. Thus
$H(u)\subseteq H(d)$ and $1=\alpha(H(u))\leq\alpha(H(d))=y$, a contradiction.
\end{proof}

\begin{remark}
No mention of the interpolation property has been made here, neither for $G$
nor for $H$.
\end{remark}

\section{A finite diagram without any lifting by interpolation groups}

The ``diagram'' of the following example is indexed not by a partially ordered set,
but by a category, namely the category with exactly one object and one nontrivial
idempotent endomorphism.

\begin{examplepf}\label{Ex:NonLiftDiag}
Let $\xt\colon\two^2\to\two^2$ defined by the rule
$\xt(\bx,\by)=(\bx\vee\by,\by)$ \pup{for all $\bx$, $\by<2$}.
Then there is no lifting $t$ of $\xt$, between interpolation groups, such
that $t^2=t$.
\end{examplepf}

\begin{proof}
Suppose otherwise. Let $G$ be an interpolation group such that $t$ is an
endomorphism of $G$. Since $G$ is an interpolation group and $G$ lifts $\two^2$,
$G$ can be decomposed as $E\times F$, where $E$ and $F$ are
\emph{simple} interpolation groups, and $t$ is given by
 \[
 t(x,y)=(\alpha x+\beta y,\gamma y),\quad
 \text{for all }(x,y)\in E\times F,
 \]
where $\alpha\colon E\to E$, $\beta\colon F\to E$, and
$\gamma\colon F\to F$ are nonzero positive homomorphisms. The
equality $t=t^2$ can then be expressed in matrix form as
 \[
 \begin{pmatrix}
 \alpha & \beta\\
 0 & \gamma
 \end{pmatrix}
=\begin{pmatrix}
 \alpha & \beta\\
 0 & \gamma
 \end{pmatrix}
\begin{pmatrix}
 \alpha & \beta\\
 0 & \gamma
 \end{pmatrix}
=\begin{pmatrix}
 \alpha^2 & \alpha\beta+\beta\gamma\\
 0 & \gamma^2
 \end{pmatrix}
 \]
that is,
 \begin{equation}\label{Eq:abcMat}
 \alpha^2=\alpha,\quad\gamma^2=\gamma,\quad\text{and}\quad
 \alpha\beta+\beta\gamma=\beta.
 \end{equation}
Hence, by repeated uses of \eqref{Eq:abcMat}, we obtain the equalities
 \begin{align*}
 \beta&=\alpha\beta+\beta\gamma\\
 &=\alpha(\alpha\beta+\beta\gamma)+(\alpha\beta+\beta\gamma)\gamma\\
 &=\alpha^2\beta+2\alpha\beta\gamma+\beta\gamma^2\\
 &=\beta+2\alpha\beta\gamma,
 \end{align*}
so we obtain that
 \begin{equation}\label{Eq:2abc=0}
 2\alpha\beta\gamma=0.
 \end{equation}
However, $E$ and $F$ are simple \poag s
and $\alpha$, $\beta$, $\gamma$ are nonzero positive homomorphisms, in
particular, $(2\alpha\beta\gamma)[F^{++}]\subseteq E^{++}$, which
contradicts \eqref{Eq:2abc=0}.
\end{proof}

\section{A diagram of finite Boolean \jzs s and
\jz-embeddings without any simplicial lifting}

\begin{examplepf}\label{Ex:NonSimplSquare}
Let $\xf\colon\two^2\hookrightarrow\two^3$ and
$\xg$, $\xh\colon\two^3\hookrightarrow\two^4$ be the maps defined by
 \begin{align*}
 \xf(\bx,\by)&=(\bx,\by,\bx\vee\by),\\
 \xg(\bu,\bv,\bw)&=(\bu,\bv,\bw,\bu\vee\bw),\\
 \xh(\bu,\bv,\bw)&=(\bu,\bv,\bw,\bv\vee\bw),
 \end{align*}
for all $\bx$, $\by$, $\bu$, $\bv$, $\bw\in\set{0,1}$. Observe that
$\xg\circ\xf=\xh\circ\xf$, see the diagram~\eqref{Eq:SemDiagr}.
 \begin{equation}\label{Eq:SemDiagr}
{
\def\labelstyle{\displaystyle}
\xymatrix{
& \two^4 & \\
\two^3\ar@{_(->}[ru]^{\xg} & & \two^3\ar@{^(->}[lu]_{\xh}\\
& \two^2\ar@{_(->}[lu]^{\xf}\ar@{^(->}[ru]_{\xf}&
}
}
 \end{equation}
Then the diagram \eqref{Eq:SemDiagr} has no lifting that uses only
simplicial \povs s, that is, a lifting of the form
 \begin{equation}\label{Eq:SimplDiagr}
{
\def\labelstyle{\displaystyle}
\xymatrix{
& \QQ^4 & \\
\QQ^3\ar[ru]^{g} & & \QQ^3\ar[lu]_{h}\\
& \QQ^2\ar[lu]^{f_0}\ar[ru]_{f_1}&
}
}
 \end{equation}
\end{examplepf}

\begin{proof}
Suppose otherwise. We identify the maps $f_0$, $f_1$, $g$, and $h$ of the
diagram \eqref{Eq:SimplDiagr} with their matrices, that have the following form:
 \begin{align*}
 f_0&=
 \begin{pmatrix}
 \alpha & 0\\ 0 & \beta\\ \xi & \eta
 \end{pmatrix},&
 f_1&=
 \begin{pmatrix}
 \alpha' & 0\\ 0 & \beta'\\ \xi' & \eta'
 \end{pmatrix},\\
 g&=
 \begin{pmatrix}
 a & 0 & 0\\ 0 & b & 0\\ 0 & 0 & c\\ u & 0 & w
 \end{pmatrix},&
 h&=
 \begin{pmatrix}
 a' & 0 & 0\\ 0 & b' & 0\\ 0 & 0 & c'\\ 0 & v' & w'
 \end{pmatrix},
 \end{align*}
with positive rational numbers $\alpha$, $\beta$, $\xi$, $\eta$,
$\alpha'$, $\beta'$, $\xi'$, $\eta'$, $a$, $b$, $c$, $u$, $w$,
$a'$, $b'$, $c'$, $v'$, $w'$. Hence the equality $g\circ f_0=h\circ f_1$
takes the form
 \[
 \begin{pmatrix}
 a & 0 & 0\\ 0 & b & 0\\ 0 & 0 & c\\ u & 0 & w
 \end{pmatrix}
 \begin{pmatrix}
 \alpha & 0\\ 0 & \beta\\ \xi & \eta
 \end{pmatrix}=
 \begin{pmatrix}
 a' & 0 & 0\\ 0 & b' & 0\\ 0 & 0 & c'\\ 0 & v' & w'
 \end{pmatrix}
 \begin{pmatrix}
 \alpha' & 0\\ 0 & \beta'\\ \xi' & \eta'
 \end{pmatrix},
 \]
that is,
 \begin{equation}\label{Eq:SimplMatr}
 \begin{pmatrix}
 a\alpha & 0\\ 0 & b\beta\\ c\xi & c\eta\\ u\alpha+w\xi & w\eta
 \end{pmatrix}=
 \begin{pmatrix}
 a'\alpha' & 0\\ 0 & b'\beta'\\ c'\xi' & c'\eta'\\
 w'\xi' & v'\beta'+w'\eta'
 \end{pmatrix}.
 \end{equation}
Hence, $\dfrac{\eta'}{\xi'}=\dfrac{c'\eta'}{c'\xi'}=
\dfrac{c\eta}{c\xi}=\dfrac{\eta}{\xi}$, while
$\dfrac{\eta'}{\xi'}=\dfrac{w'\eta'}{w'\xi'}<
\dfrac{v'\beta'+w'\eta'}{w'\xi'}=\dfrac{w\eta}{u\alpha+w\xi}<
\dfrac{w\eta}{w\xi}=\dfrac{\eta}{\xi}$, a contradiction.
\end{proof}

\begin{remark}
We have used only the last two rows of the matrices in
\eqref{Eq:SimplMatr}, the rest is there to ensure that $\xf$, $\xg$, and
$\xh$ are embeddings.

\end{remark}

\section{A semilattice map without lifting}

The following counterexample shows that Corollary~\ref{C:Nab1-1}, that states
that every \jzh\ between countable distributive \jzs s can be lifted by a
positive homomorphism between countable dimension vector spaces, cannot be
extended to uncountable semilattices:

\begin{examplepf}\label{Ex:al1diagnolift}
Let $B=\setm{x\subseteq\omega_1}%
{\text{either }x\text{ or }\omega_1\setminus x\text{ is finite}}$. So
$B$ is a Boolean semilattice of size $\aleph_1$.
Let $\xs\colon B\to\two$ be the simplest closure operator on $B$,
that is, $\xs(0)=0$ while $\xs(\bx)=1$, for all $\bx\in B\setminus\set{0}$.
Then $\xs$ has no lifting, that is, there are no positive homomorphism
$s\colon G\to S$, for some \poag s $G$ and $S$ such that $G$ has interpolation,
and no isomorphisms
$\alpha\colon B\to\Idc G$ and $\eps\colon\two\to\Idc S$ such that
the following diagram is commutative:
 \[
{
\def\labelstyle{\displaystyle}
\xymatrix{
 B\ar[r]^{\xs}\ar[d]_{\alpha} & \two\ar[d]^{\eps} \\
 \Idc G\ar[r]_{\Idc s} & \Idc S
}
}
 \]
\end{examplepf}

\begin{proof}
Suppose otherwise. In particular, $S$ is a simple \poag.
Let $\alpha(1_B)=G(u)$, for some $u\in G^{++}$. For all
$\xi<\omega_1$, the inequality $\set{\xi}< 1_B$ holds, so
$\alpha(\set{\xi})<\alpha(1_B)=G(u)$, thus there are $m_\xi\in\NN$ and
$u_\xi\in[0,m_\xi u]$ in $G$ such that $G(u_\xi)=\alpha(\set{\xi})$. There
exists an uncountable subset $U$ of $\omega_1$ such that
$m_\xi=m=\text{constant}$, for all $\xi\in U$, therefore we may replace $u$ by
$mu$ and assume that $0\leq u_\xi\leq u$, for all $\xi\in U$, thus
there exists $v_\xi\in[0,u]$ such that $u_\xi+v_\xi=u$.
Furthermore, for $\xi\neq\eta$ in $\omega_1$,
$\set{\xi}\cap\set{\eta}=\emptyset$, thus
$G(u_\xi)\cap G(u_\eta)=\set{0}$, hence $u_\xi\wedge u_\eta=0$ in $G$.

By applying this last fact together with Lemma~\ref{L:MultRef} in the
refinement monoid $G^+$ (with $I=A$ and $T_\xi=\set{0,1}$, for all $\xi\in A$) to
the system of equalities $u_\xi+v_\xi=u$ (for all $\xi\in A$), we obtain the
inequality
 \begin{equation}\label{Eq:Afin}
 \sum_{\xi\in A}u_\xi\leq u,\quad\text{for all finite }A\subset U.
 \end{equation}
Next, since $s$ lifts $\xs$, $s(u_\xi)\in S^{++}$, for all $\xi\in U$,
thus, since $S$ is simple, there exists $n_\xi\in\NN$ such that
$s(u)\leq n_\xi s(u_\xi)$. Let $V\subseteq U$ uncountable such that
$n_\xi=n=\text{constant}$, for all $\xi\in V$. Pick any distinct $\xi_0$, \dots,
$\xi_n\in V$. For all $i\in\set{0,\ldots,n}$, the inequalities
$s(u)\leq n_{\xi_i}\cdot s(u_{\xi_i})=n\cdot s(u_{\xi_i})$ hold, hence,
adding together all these inequalities, we obtain that
 \begin{align*}
 (n+1)\cdot s(u)&\leq n\cdot\sum_{i=0}^ns(u_{\xi_i})\\
 &=n\cdot s\left(\sum_{i=0}^nu_{\xi_i}\right)\\
 &\leq n\cdot s(u)\quad\text{(by \eqref{Eq:Afin})},
 \end{align*}
whence $s(u)\leq 0$, a contradiction.
\end{proof}

\section{Open problems}

We first recall the following central open problem, which is the basic motivation
of the present paper:

\begin{problem}[Problem~10.1 of \cite{GoWe}]
Let $S$ be a \jzs\ of cardinality $\aleph_1$. Does there exist a dimension group
$G$ such that $\Idc G\cong S$?
\end{problem}

A related problem is the following.

\begin{problem}\label{Pb:AnyDiag}
Can \emph{every} finite diagram (indexed by a non necessarily dismantlable
partially ordered set) of finite Boolean semilattices and
\jzh s be lifted by a diagram of dimension groups and positive homomorphisms?
\end{problem}

Example~\ref{Ex:NonLiftDiag} shows that the analogue of Problem~\ref{Pb:AnyDiag}
for diagrams indexed by arbitrary finite categories (instead of partially
ordered sets) fails. Still the problem for an arbitrary finite partially ordered
index set $P$ remains, in particular for $P=\two^3$.

We do not know about the ring-theoretical analogue of Theorem~\ref{T:dislift}:

\begin{problem}\label{Pb:IdcRings}
Let $K$ be a field. Is it the case that every diagram of finite Boolean \jzs s
that is indexed by a finite dismantlable partially ordered set can be lifted,
with respect to the $\Idc$ functor on \emph{rings}, by a diagram of locally
matricial algebras over $K$?
\end{problem}

We do not even know whether Problem~\ref{Pb:IdcRings} has a positive solution
for the square diagram of Example~\ref{Ex:NonSimplSquare}. On the other hand, it
follows from the results of F.~Weh\-rung~\cite{Wehr00} that \emph{every square
diagram of finite Boolean semilattices can be lifted, with respect to the $\Idc$
functor, by a diagram of von~Neumann regular algebras \pup{over any given
field}}.

\begin{problem}\label{Pb:Nabla1}
Let $G$ be a dimension group, let $S$ be a distributive \jzs, and let
$\xf\colon\Idc G\to S$ be a \jzh. Find sufficient conditions under which~$\xf$
can be lifted, that is, there are a dimension group $H$ and an isomorphism
$\beta\colon S\to\Idc H$ such that $\Idc f=\beta\circ\xf$?
\end{problem}

Example~\ref{Ex:nabla1fails} shows that strong restrictions on $G$, $S$, and $\xf$
are needed. Compare also with Theorem~\ref{T:Nab1-1}.

\section*{Acknowledgments}
Part of this work was completed while the second author was visiting
the department of mathematics at Charles University in Prague while staying
at the Suchdol campus. The excellent conditions provided by both places are
greatly appreciated. In particular, special thanks are due to V\'aclav
Slav\'\i k.

Another part of the work was completed while the first author visited
the University of Caen supported by the Barrande Program.
Special thanks also to Philippe Toffin.

\end{document}